\documentclass[10pt,a4paper,reqno]{amsart}
\usepackage{amsfonts,amsthm,latexsym,amsmath,amssymb,amscd,amsmath,epsf,morefloats}
\usepackage{placeins}
\usepackage{graphicx}
\usepackage{tikz}
\usetikzlibrary{arrows}

\newtheorem{theo+}           {Theorem}
\newtheorem{prop+}           {Proposition}
\newtheorem{coro+}           {Corollary}
\newtheorem{lemm+}           {Lemma}

\theoremstyle{definition}
\newtheorem{defi+}           {Definition}

\theoremstyle{remark}
\newtheorem{rema+}           {Remark}

\newenvironment{theorem}{\begin{theo+}}{\end{theo+}}
\newenvironment{proposition}{\begin{prop+}}{\end{prop+}}
\newenvironment{corollary}{\begin{coro+}}{\end{coro+}}
\newenvironment{lemma}{\begin{lemm+}}{\end{lemm+}}

\newenvironment{definition}{\begin{defi+}}{\end{defi+}}

\title{ Fingerprints, lemniscates and quadratic differentials}
\author{Alexander~Yu.~Solynin
}
\address{Department of Mathematics and Statistics, Texas Tech
University, Box 41042, Lubbock, Texas 79409}
\email{alex.solynin@ttu.edu}%

\date{\today}

\begin{document}
          \numberwithin{equation}{section}

          \begin{abstract}
We  discuss some aspects of the theory of recognition of
two-dimensional shapes by means of fingerprints of Jordan curves.
An interesting approach to problems on shape recognition suggested
by P.~Ebenfelt, D.~Khavinson, and H.~Shapiro and extended further
by M.~Younsi reveals the fact that the fingerprints of polynomial
lemniscates and, more generally, fingerprints of rational
lemniscates can be obtained as solutions to certain functional
equations involving Blaschke products.

Our main goal here is to develop an approach which relates
fingerprints of Jordan curves composed of arcs of trajectories and
orthogonal trajectories of certain quadratic differentials  with
solutions of functional equations involving pullbacks of these
quadratic differentials under appropriate Riemann mapping
functions. In particular, we show that the previous results of
P.~Ebenfelt, D.~Khavinson, and H.~Shapiro and the recent results
of M.~Younsi follow from our more general theorems as special
cases.

\end{abstract}

\maketitle

\section{Shapes and fingerprints}

\setcounter{equation}{0}

The study of two-dimensional ``shapes'' and their
``fingerprints'', initiated by A.~Kirillov \cite{K} and developed
by E.~Sharon and D.~Mumford \cite{SM}, became a rather hot topic
in recent publications on applications of complex analysis to
problems in pattern recognition.  Several authors have explored
their own ways to work in this area; see, for instance, a paper
\cite{Ma} of D.~Marshall, who used his \emph{zipper algorithm} and
a paper \cite{Wi} of  B.~Williams, who applied a \emph{circle
packing} technique. An interesting approach to the fingerprint
problem was suggested by P.~Ebenfelt, D.~Khavinson, and H.~Shapiro
in \cite{EKS}. In particular, they showed that fingerprints of
polynomial lemniscates (which, by the well-known theorem of
D.~Hilbert (cf. Chapter~4 in \cite{W}), are dense in the space of
all two dimensional shapes) are generated by solutions of
functional equations which involve Blaschke products. A simpler
proof of the main result in \cite{EKS} and its generalization to
the case of rational lemniscates were presented in a nice short
paper by M.~Younsi \cite{Y}. One more approach, which allows
certain reinterpretation of the results of P.~Ebenfelt,
D.~Khavinson, and H.~Shapiro and results of M.~Younsi, was
recently suggested by T.~Richards; see \cite{Richards} and
references therein.

My intention here is to emphasize the role of quadratic
differentials in this developing theory. In the context of image
recognition, quadratic differentials were already used by
S.~Huckemann, T.~Hotz, and A.~Munk  in their very interesting
paper \cite{HHM}. I want to stress here that our approach is
different from the approach used in \cite{HHM}. First, we
introduce, following the presentations in \cite{SM} and
\cite{EKS}, two-dimensional shapes and their fingerprints. Let
$\Gamma$ be a Jordan curve in the complex plane $\mathbb{C}$ and
let $\Omega_-$ and $\Omega_+$ denote the bounded and unbounded
components of $\overline{\mathbb{C}}\setminus \Gamma$, where
$\overline{\mathbb{C}}$ is the complex sphere. Then $\Omega_-$ and
$\Omega_+$ are simply connected domains and therefore, by the
Riemann mapping theorem, there exist conformal and one-to-one maps
$\varphi_-:\mathbb{D}\to \Omega_-$ and $\varphi_+:\mathbb{D}_+\to
\Omega_+$, where $\mathbb{D}=\{z:\,|z|<1\}$ is the unit disk and
$\mathbb{D}_+=\overline{\mathbb{C}}\setminus
\overline{\mathbb{D}}$. Figure~1 illustrates basic notations
related to the curve $\Gamma$ and mapping functions $\varphi_-$
and $\varphi_+$.


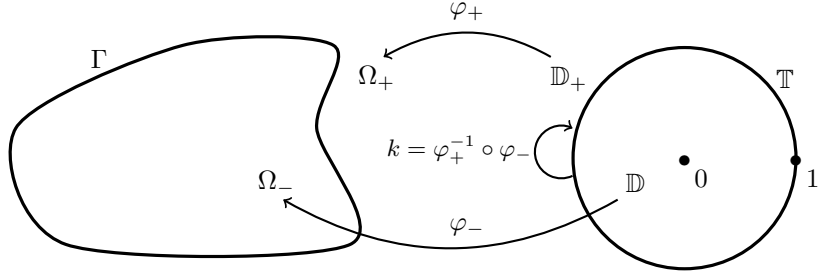
\begin{figure}
\centering %
\hspace{-2.5cm} %
\begin{minipage}{0.7\linewidth}
\begin{tikzpicture}
    [inner sep=1mm,
    place/.style={circle,draw=black!50,
    fill=blue!20,thick},
    scale=2.2]%


\draw [black, very thick](2,-0.18) circle (2/3); %
\node at (2,-0.2) [black] {$\bullet$}; %



\node at (2.1,-0.2) [below] {$0$};

\node at (8/3,-0.2) [black] {$\bullet$}; %
\node at (2.766666,-0.2) [below] {$1$};

\node at (1.6,-0.33) [right] {$\mathbb{D}$};

\coordinate (A) at (1.2,0.43);

\coordinate (B) at (0.2,0.43);

\coordinate (C) at (1.6,-0.43);

\coordinate (D) at (-0.4,-0.43);

\node at (1.15,0.3) [right] {$\mathbb{D}_+$};

\node at (2.5,0.3) [right] {$\mathbb{T}$};






\node at (-1.6,0.42) [right] {$\Gamma$};

\node at (-0.6,-0.33) [right] {$\Omega_-$};

\node at (0,0.3) [right] {$\Omega_+$};

\node at (0.7,0.6) [above] {$\varphi_+$};

\node at (0.7,-0.7) [above] {$\varphi_-$};

\node at (0.17,-0.13) [right]
{\small{$k=\varphi_+^{-1}\circ\varphi_-$}};


(1) edge[bend right] node [left] {} (2);

\path [black, thick, bend right,->]   (A) edge (B);

\path [black, thick, bend left,->]   (C) edge (D);


\draw [black, very thick] plot [smooth cycle] coordinates {(-2,0)
(-1,0.5) (-0.1,0.5) (-0.2,0) (0,-0.7) (-1.7,-0.7)};





 \draw [<-,thick] (1.33,0) arc (65:295:4.5pt);

 \end{tikzpicture}

 \medskip

\end{minipage}

\caption{Jordan curve $\Gamma$ and complementary domains
$\Omega_-$ and $\Omega_+$. 
}
\end{figure}

\FloatBarrier



We suppose that $\varphi_+$ is normalized by the conditions
$\varphi_+(\infty)=\infty$, $\varphi'_+(\infty)>0$, where
$\varphi'_+(\infty)=\lim_{z\to \infty} \varphi_+(z)/z$. The latter
normalization defines $\varphi_+$ uniquely. By the
Carath\'{e}odory theorem on  boundary correspondence (cf.
Theorem~4$'$ in \cite[Chapter~2, \S 3]{G}), each of the maps
$\varphi_-$ and $\varphi_+$ extends as a continuous one-to-one
function to the unit circle $\mathbb{T}=\partial \mathbb{D}$.
Therefore, the composition $k=\varphi_+^{-1}\circ \varphi_-$
defines an orientation preserving  automorphism of $\mathbb{T}$.
Since $\varphi_-$ is uniquely determined up to a precomposition
with a M\"{o}bius automorphism of $\mathbb{D}$, the automorphism
$k$ is also uniquely determined up to a M\"{o}bius automorphism of
$\mathbb{D}$, i.e. up to a precomposition with maps %
\begin{equation} \label{1.1} %
\phi(z)=\lambda \frac{z-a}{1-\bar{a}z}, \quad |\lambda|=1, \quad
a\in \mathbb{D}.
\end{equation} 
The equivalence class of the automorphism $k$ under the action of
the M\"{o}bius group of automorphisms (\ref{1.1}) is called the
\emph{fingerprint set} of $\Gamma$ (many authors prefer a shorter
name \emph{the fingerprint})  and elements of this set are called
\emph{fingerprints}. Furthermore, the fingerprint set of the curve
$\Gamma$  is invariant under translations and scalings of
$\Gamma$, i.e. under affine maps $L(z)=az+b$ with $a>0$, $b\in
\mathbb{C}$. The equivalence class of a Jordan curve $\Gamma$
under the action of affine maps of this form is called the
\emph{shape} and $\Gamma$ is a representative of this shape. In
what follows, we mainly work with smooth and piecewise smooth
Jordan curves $\Gamma$; i.e. we mostly work with Jordan curves
$\Gamma$,  which have a continuous or piecewise continuous unit
tangent vector. More precisely, a Jordan curve $\Gamma$ is said to
be \emph{piecewise smooth} if it consists of a finite number of
simple arcs, say  $\gamma_1,\ldots, \gamma_n$, such that every arc
$\gamma_k$ has a unit tangent vector $\vec{v}(a)$ at every point
$a\in \gamma_k$, including one-sided unit tangent vectors at the
endpoints of $\gamma_k$ and such that $\vec{v}(a)$ is continuous
on $\gamma_k$, including one-sided continuity at the endpoints of
$\gamma_k$.
 In these cases,
the corresponding shapes will be called \emph{smooth shapes} and
\emph{piecewise smooth shapes}, accordingly. Thus, we have a map
$\mathcal{F}$ from the set of all shapes into the set of all
equivalence classes of orientation preserving homeomorphisms of
$\mathbb{T}$ onto itself. Let $\mathcal{S}^1$ denote the class of
all smooth shapes in $\mathbb{C}$ and let
${\mbox{Diff}}_M(\mathbb{T})$ denote the set of equivalence
classes (under the action of the M\"{o}bius group of automorphisms
(\ref{1.1})) of orientation preserving diffeomorphisms of
$\mathbb{T}$. The following pioneering result was proved by
A.~Kirillov.

\begin{theorem}[\cite{K}]\label{Theorem-1}
The map $\mathcal{F}$ is a bijection between $\mathcal{S}^1$ and
${\mbox{Diff}}_M(\mathbb{T})$.
\end{theorem} %

In fact, Theorem~1 is a consequence of the ``fundamental theorem
of conformal welding'' proved by A.~Pfluger in 1961,
\cite{Pfluger}.

In other words, Theorem~1 says that ${\mbox{Diff}}_M(\mathbb{T})$
parameterizes the set $\mathcal{S}^1$ of all smooth shapes. The
set of all diffeomorphisms of $\mathbb{T}$ is rather large. So,
P.~Ebenfelt, D.~Khavinson, and H.~Shapiro \cite{EKS} studied
possible parameterizations of fingerprints of polynomial
lemniscates.  A lemniscate of a polynomial $P(z)$ at level $c>0$
is defined as $L_P(c)=\{z:\,|P(z)|=c\}$. Although in this paper we
deal mostly  with analytic and connected lemniscates, we emphasize
here that the polynomial lemniscates are not necessarily connected
or Jordan, in general. Later on, we will consider also rational
lemniscates and lemniscates associated with nonconstant
meromorphic functions. The set of analytic and connected
polynomial lemniscates is much smaller than the set of all smooth
Jordan curves, nevertheless, by the Hilbert theorem mentioned
above, it still can be used to approximate any Jordan shape. The
following theorem of P.~Ebenfelt, D.~Khavinson, and H.~Shapiro
characterizes exactly which elements of the set of diffeomorphisms
of the unit circle  $\mathbb{T}$   appear to be the fingerprints
of
polynomial lemniscates. %

\begin{theorem}[\cite{EKS}] \label{Theorem-2}
Let $P(z)=c_n z^n+c_{n-1}z^{n-1}+\ldots +c_0$ be a polynomial of
degree $n$ with $c_n>0$ such that $L_P(1)$ is analytic and
connected and let $k: \mathbb{T}\to \mathbb{T}$ be a fingerprint
of $L_P(1)$. Then  $k(z)$  is given by the equation %
\begin{equation} \label{1.2} %
k(z)={B(z)}^{1/n},
\end{equation} %
where $B(z)$ is a Blaschke product of degree $n$, %
$$ 
B(z)=e^{i\alpha}\prod_{k=1}^n \frac{z-a_k}{1-\overline{a_k}z}, %
$$ 
with some real $\alpha$, where $a_k=\varphi_-^{-1}(\zeta_k)$ and
$\zeta_1,\ldots,\zeta_n$ are the zeroes of $P(z)$ counting
multiplicities.

Conversely, given any Blaschke product of degree $n$, there is a
polynomial $P(z)$ of the same degree whose lemniscate $L_P(1)$ is
analytic and connected and has $k(z)={B(z)}^{1/n}$ as its
fingerprint. Moreover, $P(z)$ is unique up to precomposition with
an affine map of the form $L(z)=az+b$ with $a>0$ and $b\in
\mathbb{C}$.
\end{theorem} %

Throughout the paper, the functional notations like $k=k(z)$,
$k=k(e^{i\theta})$, etc. will be reserved to denote  fingerprints.
However, the letter $k$ will be also used as an index in several
formulas involving products and sums. We believe that this usage
will not confuse the reader since the meaning of symbol $k$ will
be clear from the context.

The proof of Theorem~2 given in \cite{EKS} is rather involved. A
shorter proof was given by M.~Younsi \cite{Y} who also proved a
counterpart of this theorem for the case of rational lemniscates,
which was conjectured in \cite{EKS}.

\begin{theorem}[\cite{Y}] \label{Theorem-3}
Let $R(z)$ be a rational function of degree $n$ with
$R(\infty)=\infty$ such that its lemniscate $L_R(1)=\{z:\,
|R(z)|=1\}$ is analytic and connected and  let $k: \mathbb{T}\to
\mathbb{T}$ be a fingerprint of $L_R(1)$. Then  $k(z)$
is given by a solution to the functional equation%
\begin{equation} \label{1.4} %
A\circ k=B,
\end{equation} %
where $A(z)$ and $B(z)$ are Blaschke products of degree $n$ and
$A(\infty)=\infty$.%

Conversely, given any solution $k(z)$ to a functional equation
$A\circ k=B$, where $A(z)$ and $B(z)$ are Blaschke products of
degree $n$ and $A(\infty)=\infty$, there exist a rational function
$R(z)$ of  degree $n$ with $R(\infty)=\infty$ whose lemniscate
$L_R(1)$ is analytic and connected and has $k(z)$ as its
fingerprint. %
\end{theorem} %

We will see later that equations (\ref{1.2}) and (\ref{1.4}) can
be interpreted in a rather natural way in terms of quadratic
differentials. Necessary properties of  quadratic differentials
and their relation with fingerprints are discussed in Section~2.
In particular, in Lemma~1 and Theorem~4 we show how to construct a
fingerprint of a piecewise smooth Jordan curve consisting of arcs
of trajectories and/or arcs of orthogonal trajectories of certain
quadratic differentials. Then in Theorems~5 and 6 we describe the
so-called \emph{welding} procedure used to construct  a curve by
its fingerprint related to a pair of quadratic differentials, one
of which is defined in the unit disk and the other one is defined
in the exterior of the unit disk.

In fact, the welding procedure mentioned above can be used to
construct a rather broad variety of quadratic differentials with
desired topological properties.  This approach will be explained
in Section~3. Our presentation of the results in Section~3 follows
closely to the exposition of results in our paper \cite{S1} where
this approach was initiated. Our main goal in this section is to
set up terminology and present results in the form which might be
useful in future work on image recognition.

In Section~4, we work with lemniscates of meromorphic functions.
In particular, we show  how earlier results of P.~Ebenfelt,
D.~Khavinson, and H.~Shapiro in \cite{EKS} and more recent results
of M.~Younsi in \cite{Y} follow from our theorems presented in
Section~2.

Finally in Section~5, we will demonstrate how our method can be
applied to study fingerprints of polygonal Jordan curves when the
corresponding Riemann mapping functions can be expressed in terms
of the well-known Schwarz-Christoffel integrals.

\medskip

Some interesting problems related to our results presented in this
paper are discussed in the  manuscript \cite{FKV} by A.~Frolova,
D.~Khavinson, and A.~Vasil'ev, which was communicated to us by one
of the referees.

\medskip

\textbf{Acknowledgements.} The author would like to thank
anonymous referees for careful reading of this paper and several
valuable suggestions, which, in our opinion,  resulted in a more
clear presentation of some portions of this paper. Also the author
would like to thank the referees for bringing papers
\cite{Pfluger} and \cite{Richards} and manuscript \cite{FKV} to
his attention.

\section{Fingerprints and  quadratic differentials}

\setcounter{equation}{0}

For the definitions and results on quadratic differentials needed
for
 the purposes of this paper the reader may consult the classical monographs of J.A.~Jenkins \cite{J2}
and K.~Strebel \cite{St} and  papers \cite{J1}, \cite{S1}, and
\cite{S2}.

A quadratic differential on a domain $D\subset
\overline{\mathbb{C}}$ is a differential form $Q(z)\,dz^2$ with
meromorphic $Q(z)$ and with the conformal transformation rule
\begin{equation} \label{2.1} %
Q_1(\zeta)\,d\zeta^2=Q(\varphi(z))\left(\varphi'(z)\right)^2\,dz^2,
\end{equation}
where $\zeta=\varphi(z)$ is a conformal map from $D$ onto a domain
$G$ in the extended plane of the parameter $\zeta$. Then the zeros
and poles of $Q(z)$ are critical points of $Q(z)\,dz^2$, in
particular, zeros and simple poles are finite critical points and
poles of order greater than $1$ are infinite critical points of
$Q(z)\, dz^2$. A trajectory (respectively, orthogonal trajectory)
of $Q(z)\,dz^2$ is a closed analytic Jordan curve or maximal open
analytic arc $\gamma\subset D$ such
that %
$$ 
Q(z)\,dz^2>0 \quad {\mbox{along $\gamma$}} \quad \quad
({\mbox{respectively}}, Q(z)\,dz^2<0 \quad {\mbox{along
$\gamma$}}).
$$ 

A trajectory $\gamma$ is called \emph{critical} if at least one of
its endpoints is a finite critical point of $Q(z)\,dz^2$. If
$\gamma$ is a rectifiable arc in $D$ then its $Q$-length is
defined by $|\gamma|_Q=\int_\gamma |Q(z)|^{1/2}\,|dz|$.  An
important property of quadratic differentials is that the
transformation rule (\ref{2.1}) respects trajectories, orthogonal
trajectories and their $Q$-lengths, as well as it respects
critical points together with their multiplicities and trajectory
structure nearby.

\medskip

Returning to fingerprints, suppose that $\Gamma$ is a piecewise
smooth Jordan curve in the plane $\mathbb{C}$ with the parameter
$\zeta$ and that $Q(\zeta)\,d\zeta^2$ is a quadratic differential
on some neighborhood $G$ of $\Gamma$. The best case scenario is
when $G=\overline{\mathbb{C}}$. More generally, we may assume,
without loss of generality, that $G$ is a doubly connected domain
bounded by Jordan analytic  curves and that $\Gamma$ separates
boundary components of $G$. Let $Q_-(z)\,dz^2$ and $Q_+(z)\,dz^2$
denote pullbacks of $Q(\zeta)\,d\zeta^2$ under the conformal maps
$\zeta=\varphi_-(z)$ and $\zeta=\varphi_+(z)$ defined in Section~1
and let $\tau_-(\zeta)=\varphi_-^{-1}(\zeta)$ and
$\tau_+(\zeta)=\varphi_+^{-1}(\zeta)$.
Then  %
\begin{equation} \label{2.3} %
Q(\zeta)\,d\zeta^2=Q_-(\tau_-(\zeta))(\tau'_-(\zeta))^2\,d\zeta^2
\quad
{\mbox{for $\zeta\in (\Omega_-\cap G)\cup \Gamma$}} %
\end{equation} %
and %
\begin{equation} \label{2.4} %
 Q(\zeta)\,d\zeta^2=Q_+(\tau_+(\zeta))(\tau'_+(\zeta))^2\,d\zeta^2 \quad
{\mbox{for $\zeta\in (\Omega_+\cap G)\cup \Gamma$}}.
\end{equation} %

Since $\Gamma$ is piecewise smooth it follows that  each of the
maps $\tau_-(\zeta)$ and $\tau_+(\zeta)$ and their derivatives
$\tau'_-(\zeta)$ and $\tau'_+(\zeta)$ can be extended by
continuity to  any smooth arc of $\Gamma$.  If $\Gamma$ is smooth
at $\zeta$  then equations (\ref{2.3}) and
(\ref{2.4}) imply that %
\begin{equation} \label{2.4.1} %
Q_-(\tau_-(\zeta))(\tau'_-(\zeta))^2\,d\zeta^2=Q_+(\tau_+(\zeta))(\tau'_+(\zeta))^2\,d\zeta^2.
\end{equation} %
Changing the variable in (\ref{2.4.1}) via $\zeta=\varphi_-(z)$,
we
obtain  the equivalent equation%
\begin{equation} \label{2.4.2} %
Q_-(z)\,dz^2=Q_+(k(z))(k'(z))^2\,dz^2,
\end{equation} %
which holds for all points $z\in \mathbb{T}$ such that
$\varphi_-(z)$ belongs to a smooth arc of $\Gamma$. Here
$k=\varphi_+^{-1}\circ \varphi_-$ is  an homeomorphism from
$\mathbb{T}$ onto itself and therefore it is a  fingerprint of
$\Gamma$ as we defined  in Section~1.

Taking square roots of both sides of
equation~(\ref{2.4.2}) and then integrating   along  the unit circle, we  obtain %
\begin{equation} \label{2.6} %
\int_{z_0}^z \sqrt{Q_+(k(\tau))}k'(\tau)\,d\tau= \int_{z_0}^z \sqrt{Q_-(\tau)}\,d\tau, %
\end{equation} %
where %
\begin{equation} \label{2.6.1} %
z_0=e^{i\theta_0}, \quad  z=e^{i\theta} \quad {\mbox{with $0\le
\theta_0< 2\pi$, \ \ $\theta_0\le \theta\le \theta_0+2\pi$}} %
\end{equation} %
and integration in (\ref{2.6})   is taken along a circular arc
$l(\theta_0,\theta)=\{\tau=e^{it}:\,\theta_0\le t\le \theta\}$.

\medskip


We note that equation (\ref{2.6}) contains integrals of square
roots of meromorphic functions and therefore it requires
appropriate interpretation.  In the case when $\Gamma$ does not
contain infinite critical points of  $Q(\zeta)\,d\zeta^2$, the
integrals in (\ref{2.6}) are finite and therefore this equation
holds for all $z_0$ and $z$ as in (\ref{2.6.1}) if the branches of
square roots in (\ref{2.6}) are chosen appropriately.

If $\Gamma$ contains infinite critical points of
$Q(\zeta)\,d\zeta^2$ then (\ref{2.6}) holds  if the arc of
integration $l(\theta_0,\theta)$ does not contain infinite
critical points of $Q_-(z)\,dz^2$. In the presence of infinite
critical points of $Q_-(z)\,dz^2$ on $\mathbb{T}$, the integrals
in (\ref{2.6}) blow up and we actually have to consider  separate
equations for each maximal arc of $\mathbb{T}$ which is free of
infinite critical points.

\medskip

The following lemma summarizes  the simple observations made above
for the case when infinite critical points are absent.

\begin{lemma} \label{Lemma-1} %
Let $\Gamma$ be a piecewise smooth Jordan curve and let
$Q_-(z)\,dz^2$ and $Q_+(z)\,dz^2$ be pullbacks of the quadratic
differential $Q(\zeta)\,d\zeta^2$ introduced above. If $\Gamma$
does not contain infinite critical points of $Q(\zeta)\,d\zeta^2$
then there is a  fingerprint $k:\mathbb{T}\to \mathbb{T}$ of
$\Gamma$  given by a solution to the
functional equation %
\begin{equation} \label{2.4.3} %
\mathcal{A}\circ k=\mathcal{B}, %
\end{equation} %
where %
\begin{equation} \label{2.4.5} %
\mathcal{A}(z)=\int_{k(z_0)}^z \sqrt{Q_+(\tau)}\,d\tau, \quad
\quad \mathcal{B}(z)=\int_{z_0}^z
\sqrt{Q_-(\tau)}\,d\tau %
\end{equation}
with  $z_0$ as in (\ref{2.6.1}) and with appropriately chosen
branches of the radicals.
\end{lemma}

Since solutions to equations (\ref{2.4.2}), (\ref{2.6}),   and
(\ref{2.4.3}) provide a kind of welding of $\Gamma$ by quadratic
differentials,  we will call them  equations of \emph{QD-welding}
of $\Gamma$. These equations do not provide much information on
the fingerprint of $\Gamma$ unless we know how to construct
quadratic differentials $Q_-(z)\,dz^2$ and $Q_+(z)\,dz^2$ with the
required properties. Fortunately enough, there is a rather general
procedure, suggested in our paper \cite{S2}, which in many
important cases allows us to construct the quadratic differentials
with the desired properties. This construction, which  presents
some independent interest, will be discussed in Section~3.

Equations (\ref{2.4.3}) and (\ref{1.4}) have the same form and
this, of course, is not a coincidence. In Section~4, we will show
that equations  (\ref{1.2}) and (\ref{1.4})
 are special forms of equation  (\ref{2.4.3}) for the cases of
polynomial and rational lemniscates, respectively.

In Lemma~1 we assumed that $\Gamma$ does not contain infinite
critical points of  $Q(\zeta)\,d\zeta^2$. To address the
fingerprint problem in the general case, we need an approximation
lemma, which easily follows from the classical convergence theorem
of T.~Rad\'{o}. We need the following definition, which is
included in the statement of Theorem~2 in \cite[Chapter~2,
\S5]{G}.

\begin{definition} %
We say that a sequence of Jordan curves $\Gamma_n$,
$n=1,2,\ldots$,  converges to a Jordan curve $\Gamma$ in
Rad\'{o}'s sense if for every $\varepsilon>0$ there is positive
integer $N$ such that for every $n\ge N$ there exists a continuous
one-to-one correspondence between the points of $\Gamma_n$ and the
points of $\Gamma$ such that the distance between corresponding
points of $\Gamma_n$ and $\Gamma$ is less than $\varepsilon$.
\end{definition}

\begin{lemma} \label{Lemma-1.1} %
Let $\Gamma_n$, $n=1,2,\ldots$, be a sequence of Jordan curves
which converges in Rad\'{o}'s sense to a  Jordan curve $\Gamma$.
Let $\Omega_-^n$ and $\Omega_+^n$ denote the bounded and unbounded
components of $\overline{\mathbb{C}}\setminus \Gamma_n$ and let
$z_0\in \Omega_-$ and $z_0\in \Omega_-^n$ for all $n=1,2,\dots$
Let $k_n=\left(\varphi_+^n\right)^{-1}\circ \varphi_-^n$  and
$k=\varphi_+^{-1}\circ \varphi_-$  denote the fingerprints of
$\Gamma_n$ and $\Gamma$ unequally determined by conditions
$\varphi_-^n(0)=z_0$, $\left(\varphi_-^n\right)'(0)>0$ and
$\varphi_-(0)=z_0$, $\left(\varphi_-\right)'(0)>0$, respectively.
Then, $k_n(e^{i\theta})$ converges to $k(e^{i\theta})$ uniformly
on $\mathbb{T}$.
\end{lemma} %

\noindent %
\emph{Proof. } 
Rad\'{o}'s convergence theorem (see Theorem~2 in \cite[Chapter~2,
\S5]{G}) implies that $\varphi_-^{n}(z)\to \varphi_-(z)$ uniformly
on $\overline{\mathbb{D}}$ and that $\varphi_+^{n}(z)\to
\varphi_+(z)$ uniformly on the closed annulus $\{z:\,1\le |z|\le
2\}$. Since $k_n=\left(\varphi_+^{n}\right)^{-1}\circ
\varphi_-^{n}$ and $k=\varphi_+^{-1}\circ \varphi_-$, the lemma
follows. \hfill $\Box$

Now, we can extend Lemma~1 for the case when a piecewise smooth
Jordan curve $\Gamma$ contains infinite critical points of
$Q(\zeta)\,d\zeta^2$ as follows. For any $a\in \Gamma$ and
$\varepsilon>0$ small enough, let $\gamma_\varepsilon(a)$ denote
an arc of  $\Gamma$ lying inside the circle
$C_\varepsilon(a)=\{\zeta:\,|\zeta-a|=\varepsilon\}$. Since
$\Gamma$ is piecewise smooth the arc $\gamma_\varepsilon(a)$ is
defined uniquely for all $\varepsilon>0$ small enough.
Furthermore,  let $l_\varepsilon(a)$ denote an arc of
$C_\varepsilon(a)$, which lies in the closure of the domain
$\Omega_-$ and joins the endpoints of $\gamma_\varepsilon(a)$.
Since $\Gamma$ is piecewise smooth the arc $l_\varepsilon(a)$ is
also defined uniquely for all $\varepsilon>0$ small enough.

Let $\Gamma_\varepsilon$ denote the curve obtained from $\Gamma$
by replacing $\gamma_\varepsilon(a)$ with $l_\varepsilon(a)$ for
all infinite critical points $a\in \Gamma$. If $\Gamma$ is
piecewise smooth and $\varepsilon>0$ is small enough, then
$\Gamma_\varepsilon$ is a piecewise smooth Jordan curve which does
not contain infinite critical points of $Q(\zeta)\,d\zeta^2$. Take
a sequence $\varepsilon_n\to 0$. By Lemma~1, the fingerprint $k_n$
of $\Gamma_{\varepsilon_n}$ can be obtained from
equation~(\ref{2.4.3}). Then, by Lemma~2, the fingerprint $k$ of
$\Gamma$ can be found as the limit %
\begin{equation} \label{corrections p6}
k(e^{i\theta})=\lim_{n\to \infty} k_n(e^{i\theta}). %
\end{equation} %

We note here that Jordan curves $\Gamma_\varepsilon$ described
above can be constructed in a more general case when $\Gamma$ is
smooth except at most finite number of points $a\in \Gamma$ such
that for each of these points and every $\varepsilon>0$
sufficiently small each of the intersections $\Omega_-\cap
C_\varepsilon(a)$ and $\Omega_+\cap C_\varepsilon(a)$ consists of
a single arc.  In particular, this happens in the case when
$\Gamma$ consists of a finite number of arcs of trajectories
and/or orthogonal trajectories of a quadratic differential. We
emphasize that, in the latter case, $\Gamma$ may have points $a\in
\Gamma$ such that in a sufficiently small neighborhood of $a$ the
arcs of $\Gamma\setminus \{a\}$ behave like logarithmic spirals.
Therefore, the limit relation (\ref{corrections p6}) remains valid
if $\Gamma$ is not necessarily piecewise smooth but consists of a
finite number of arcs of trajectories and/or orthogonal
trajectories of some quadratic differential (as in our Theorems~4
and 5 presented below).

\medskip

We admit here that equation~(\ref{2.4.3}) is of little practical
use unless we know how to find functions $\mathcal{A}(z)$ and
$\mathcal{B}(z)$ defined by (\ref{2.4.5}) or, equivalently, how to
construct quadratic differentials $Q_-(z)\,dz^2$ and
$Q_+(z)\,dz^2$.  In general, such construction looks problematic.
Rare cases when this is possible include cases of polynomial and
rational lemniscates presented in Theorems \ref{Theorem-2} and
\ref{Theorem-3}. Below in this section, we explain how to find a
general form of the quadratic differentials $Q_-(z)\,dz^2$ and
$Q_+(z)\,dz^2$ making some additional assumptions.

Namely, we suppose now that a quadratic differential
$Q(\zeta)\,d\zeta^2$ is defined on the whole complex sphere
$\overline{\mathbb{C}}$. Then, of course, $Q(\zeta)$ is a rational
function. Furthermore, we suppose that $\Gamma$ is a Jordan curve
consisting of a finite number of arcs, $\gamma_1$, \ldots,
$\gamma_m$,  of trajectories and/or orthogonal trajectories of
$Q(\zeta)\,d\zeta^2$ and their endpoints. In this context, the
whole trajectory or orthogonal trajectory is also considered as an
arc. In particular, we allow cases when $\Gamma$ consists of a
single trajectory or orthogonal trajectory, closed or not.
Infinite critical points on $\Gamma$ are also allowed.

First, we discuss briefly the possible structure of $\Gamma$ in a
neighborhood of a point $\zeta_0\in \Gamma$ where two arcs, say
$\gamma_1$ and $\gamma_2$, meet. All properties stated in items
(1)--(5) below follow from the well-known results on the local
structure of trajectories of quadratic differentials; see, for
instance, Ch.~3 in \cite{J2}.%
\begin{enumerate} %
\item[(1)] %
 If $\zeta_0$ is a regular point of
$Q(\zeta)\,d\zeta^2$, then one of the arcs $\gamma_1$ or
$\gamma_2$ is an arc of a trajectory and the other one is an arc
of an orthogonal trajectory. In this case, $\gamma_1$ and
$\gamma_2$ form a corner of opening $\pi/2$ with respect to one of
the domains $\Omega_-$ and $\Omega_+$ and a corner of opening
$3\pi/2$ with respect to the other one. %
\item[(2)] %
 If $\zeta_0$ is a zero of $Q(\zeta)\,d\zeta^2$ of order $n$ then
the arcs $\gamma_1$ and $\gamma_2$ form an angle of opening $\pi
k/(n+2)$ with some integer $k$, $0<k<2(n+2)$. Moreover, if $k$ is
odd then one of the arcs $\gamma_1$ or $\gamma_2$ is an arc of a
trajectory and the other one is an arc of an orthogonal
trajectory. If $k$ is even then  both $\gamma_1$ and $\gamma_2$
are either arcs of trajectories or arcs of orthogonal
trajectories. %
\item[(3)] %
If $\zeta_0$ is a simple pole, then one of the arcs $\gamma_1$ or
$\gamma_2$ is an arc of a trajectory and the other one is an arc
of an orthogonal trajectory and $\Gamma$ is analytic at $\zeta_0$. %
\item[(4)] %
If  $\zeta_0$ is a pole of order two then both  $\gamma_1$ and
$\gamma_2$ are  either arcs of trajectories or arcs of orthogonal
trajectories. If they are arcs of trajectories then
$Q(\zeta)\,dz^2$ has a radial or spiral structure of trajectories
near $\zeta_0$ and if they are arcs of orthogonal trajectories
$Q(\zeta)\,dz^2$ has circular or spiral structure of trajectories.
In case of the radial or circular structure of trajectories,
$\gamma_1$ and $\gamma_2$ can form any angle $\alpha$,
$0<\alpha<2\pi$. In case of the spiral structure of trajectories,
$\gamma_1$ and $\gamma_2$ in a neighborhood of $\zeta_0$ look like
logarithmic spirals and do not form any angle at $\zeta_0$. %
\item[(5)] %
 Finally, if $\zeta_0$ is a pole of order $n\ge 3$, then
$\gamma_1$ and $\gamma_2$ can form any angle of opening $\pi
k/(n-2)$ with some integer $k$, $0\le k\le 2(n-2)$. Moreover, if
$k$ is odd then one of the arcs $\gamma_1$, or $\gamma_2$ is an
arc of a trajectory and the other one is an arc of an orthogonal
trajectory. If $k$ is even then  both $\gamma_1$ and $\gamma_2$
are either arcs of trajectories or arcs of orthogonal
trajectories. %
\end{enumerate} %

Summarizing, we conclude that $\Gamma$ may contain arcs of
logarithmic spirals, and therefore it is not necessarily
rectifiable, and it may have interior and exterior cusps and
corners of any angle.



 %

%

\medskip
\FloatBarrier

Now, we discuss pullbacks $Q_-(z)\,dz^2$ and $Q_+(z)\,dz^2$ of the
quadratic differential $Q(\zeta)\,d\zeta^2$. It is enough to
consider the quadratic differential $Q_-(z)\,dz^2$, the trajectory
structure of $Q_+(z)\,dz^2$ is similar. Under our assumptions,
$Q_-(z)\,dz^2$ is defined on the whole unit disk $\mathbb{D}$.
Furthermore, the unit circle $\mathbb{T}$ consists of a finite
number of arcs of trajectories and/or orthogonal trajectories of
$Q_-(z)\,dz^2$ and their endpoints. Since $Q_-(z)\,dz^2$ is real
on $\mathbb{T}$ (possibly except a finite number of points), it
follows from the reflection principle that $Q_-(z)\,dz^2$ can be
continued to a quadratic differential on $\overline{\mathbb{C}}$,
for which we will keep the same notation $Q_-(z)\,dz^2$, and its
trajectory structure is symmetric with respect to $\mathbb{T}$. In
particular, the set of zeros and the set of poles of $Q_-(z)\,
dz^2$ are both symmetric with respect to $\mathbb{T}$.

Since conformal mappings preserve critical points of quadratic
differentials as well as the local structure of their
trajectories, it follows that the set of critical points of
$Q_-(z)\,dz^2$ inside $\mathbb{D}$ is in a one-to-one
correspondence with the set of critical points of
$Q(\zeta)\,d\zeta^2$ inside $\Omega_-$ and that $Q_-(z)\,dz^2$ and
$Q(\zeta)\,d\zeta^2$ have the same local structure near
corresponding critical points.

Now, we identify possible critical points of $Q_-(z)\,dz^2$ on
$\mathbb{T}$. All our conclusions follow from our knowledge of the
 local structure of trajectories of $Q(\zeta)\,d\zeta^2$ near endpoints of arcs composing   $\Gamma$
 (and therefore from our knowledge of local structure of trajectories of $Q_-(z)\,dz^2$ near corresponding points of
 $\mathbb{T}$). When working with critical points of
 $Q_-(z)\,dz^2$ situated on $\mathbb{T}$, we consider $\Gamma$ as
 the boundary of $\Omega_-$ oriented counterclockwise.  To
 emphasize this fact, we use notation $\Gamma_-$ instead of
 $\Gamma$ in items (a)--(g) below. Also, all angles mentioned in
 items (a)--(g) are considered with respect to the domain
 $\Omega_-$.

 \begin{enumerate} %
 \item[(a)] %
 If $\zeta_0$ is a regular point of $Q(\zeta)\,d\zeta^2$ where
 $\Gamma_-$  forms an angle $\pi/2$, then $z_0=\tau_-(\zeta_0)$
 is a simple pole of $Q_-(z)\,dz^2$.
 \item[(b)] %
 If $\zeta_0$ is a regular point of $Q(\zeta)\,d\zeta^2$ where
 $\Gamma_-$  forms an angle $3\pi/2$, then $z_0=\tau_-(\zeta_0)$
 is a simple zero of $Q_-(z)\,dz^2$.
 \item[(c)] %
 If $\zeta_0$ is a zero of $Q(\zeta)\,d\zeta^2$ of order $n$ where
 $\Gamma_-$  forms an angle $\frac{\pi k}{n+2}$, then $z_0=\tau_-(\zeta_0)$
 is a simple pole of $Q_-(z)\,dz^2$ if $k=1$ , a regular point if $k=2$, and a zero of order $k-2$ if $3\le k\le 2n+3$.
 \item[(d)] %
  If $\zeta_0$ is a simple pole of $Q(\zeta)\,dz^2$, then $z_0=\varphi_-(\zeta_0)$
 is a simple pole of $Q_-(z)\,dz^2$.
 \item[(e)] %
 If $\zeta_0$ is a pole of order $2$ of $Q(\zeta)\,dz^2$ and $\gamma_1$ and $\gamma_2$ are arcs of orthogonal trajectories,
  then $z_0=\varphi_-(\zeta_0)$
 is a pole  of $Q_-(z)\,dz^2$ of order $2$ with  circular trajectory
 structure.
 \item[(f)] %
 If $\zeta_0$ is a pole of order $2$ of $Q(\zeta)\,dz^2$ and $\gamma_1$ and  $\gamma_2$ are arcs of trajectories,
  then $z_0=\varphi_-(\zeta_0)$
 is a pole  of $Q_-(z)\,dz^2$ of order $2$ with a radial trajectory structure.
\item[(g)] %
 If $\zeta_0$ is a pole of order  $n\ge 3$ of $Q(\zeta)\,dz^2$ where
 $\Gamma_-$  forms an angle $\frac{\pi k}{n-2}$,  then $z_0=\varphi_-(\zeta_0)$
 is a pole  of $Q_-(z)\,dz^2$ of order $k+2$, $0\le k\le 2(n-2)$. In particular,
 if $k=0$ and $\gamma_1$ and  $\gamma_2$ are arcs of trajectories,
 then $z_0$ is a pole of $Q_-(z)\,dz^2$ of order $2$ with radial trajectory structure
 and if $k=0$ and $\gamma_1$, $\gamma_2$ are arcs of orthogonal trajectories,
 then $z_0$ is a pole of $Q_-(z)\,dz^2$ of order $2$ with the circular trajectory structure.
\end{enumerate} %

All statements (a)--(g) follow from the well-known results on the
local structure of trajectories of quadratic differentials. To
demonstrate details, we consider, for instance, item (b). In all
other cases the argument is similar and therefore is left to the
interested reader. In case (b), there is a neighborhood $G\subset
\Omega_-$ of $\zeta_0$ (considered as a boundary point of
$\Omega_-$), which is bounded by three arcs, say $l_1$, $l_2$, and
$l_3$, of trajectories and  three arcs, say $s_1$, $s_2$, and
$s_3$, of orthogonal trajectories of $Q(\zeta)\,d\zeta^2$ as it is
shown in Figure~2. This neighborhood is filled up with arcs of
trajectories of $Q(\zeta)\,d\zeta^2$, one of which, let $l'_1$, is
a continuation of a boundary arc $l_1$ through the point
$\zeta_0$. Let $G_z^-$ be the image of $G$ under the mapping
$\tau_-=\varphi_-^{-1}$ and let $G_z$ denote the extension of the
domain $G_z^-$ obtained by reflection with respect to the arc
$L_1\cup S_1\cup \{z_0\}$  of the unit circle. Here
$L_1=\tau_-(l_1)$, $S_1=\tau_-(s_1)$, $z_0=\tau_-(\zeta_0)$. Then
$z_0$ is a point of $G_z$, which serves as an end point for three
critical trajectories, $L_1$, $L'_1=\tau_-(l'_1)$, and
$L''_1=\{z:\,1/\bar z \in L'_1\}$,  of the quadratic differential
$Q_-(z)\,dz^2$ as it is shown in Figure~2. Since no other
trajectory of $Q_-(z)\,dz^2$ terminates at $\tau_-(\zeta_0)$ it
follows from the known results on the local structure of
trajectories (see, for instance,  Theorem~3.2 in \cite{J2}) that
$z_0=\tau_-(\zeta_0)$ is a simple zero of $Q_-(z)\,dz^2$.

 As we have mentioned above the nature
of critical points of $Q_+(z)\,dz^2$ on $\mathbb{T}$ is similar.



 %



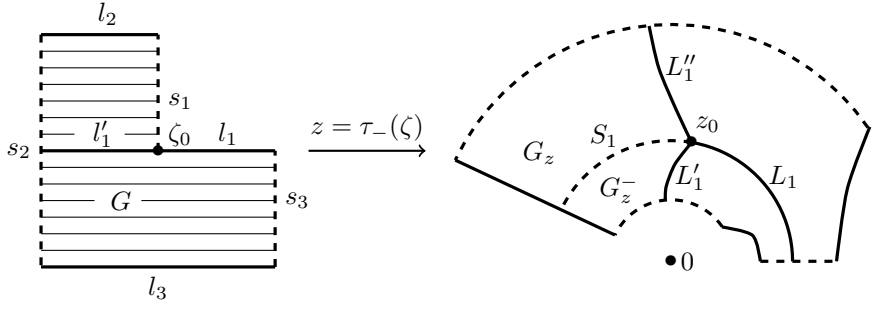
\begin{figure}
\hspace{0.2cm} %
\begin{minipage}{0.10\linewidth}
\begin{tikzpicture} [thick,scale=0.22] 
\draw [black, very thick] (0,0) --(7,0); %

\draw [black, very thick] (0,7)--(-7,7); %
\draw [black, thin] (0,6)--(-7,6); %
\draw [black, thin] (0,5)--(-7,5); %
\draw [black, thin] (0,4)--(-7,4); %
\draw [black, thin] (0,3)--(-7,3); %
\draw [black, thin] (0,2)--(-7,2); %
\draw [black, thin] (0,1)--(-2,1); %
\draw [black, thin] (-7,1)--(-5,1); %
\draw [black, very thick] (0,0)--(-7,0); %

\draw [black, very thick] (-7,-7)--(7,-7); %
\draw [black, thin] (-7,-6)--(7,-6); %
\draw [black, thin] (-7,-5)--(7,-5); %
\draw [black, thin] (-7,-4)--(7,-4); %
\draw [black, thin] (-7,-3)--(-3.8,-3); %
\draw [black, thin] (-1,-3)--(7,-3); %
\draw [black, thin] (-7,-2)--(7,-2); %
\draw [black, thin] (-7,-1)--(7,-1); %
\draw [black, dashed,very thick](0,0) --(0,7); %
\draw [black,dashed,very thick] (-7,7) --(-7,-7); %
\draw [black,dashed,very thick] (7,-7) --(7,0); %

\node at (0,0)  {$\bullet$}; %


\node at (-1,-3) [left] {$G$};  %
\node at (0,1) [right] {$\zeta_0$}; %
 \node at (3,1) [right] {$l_1$}; %
  \node at (-4.5,1) [right] {$l'_1$}; %
 \node at (-3,7) [above] {$l_2$}; %
  \node at (0,-7) [below] {$l_3$}; %

  \node at (0,3) [right] {$s_1$}; %
   \node at (-7,0) [left] {$s_2$}; %
    \node at (7,-3) [right] {$s_3$}; %

\draw [->, thick] (9,0)--(16,0); %

\node at (12.5,0) [above] {$z=\tau_-(\zeta)$}; %

\end{tikzpicture} %
\end{minipage} %
\hspace{4.5cm} %
\begin{minipage}{.5\textwidth}
\begin{tikzpicture} [thick,scale=1.6] , 




\draw [black, very thick]
({e^(0.674394308508077)*cos(2.7*180/pi)},
{e^(0.674394308508077)*sin(2.7*180/pi)})--({e^(-0.674394308508077)*cos(2.7*180/pi)},
{e^(-0.674394308508077)*sin(2.7*180/pi)}); %

\draw [black, dashed, very thick] (0.7346390419363942,
0)--(1/0.7346390419363942,0);

 \draw [black,very thick,domain=0:1.4*180/pi] plot ({cos(\x)}, {sin(\x)}); %

 \draw [black,dashed, very thick,domain=1.4*180/pi:2.7*180/pi] plot ({cos(\x)}, {sin(\x)}); %

 \node at ({cos(1.4*180/pi)}, {sin(1.4*180/pi)})  {$\bullet$}; %

  \draw [black,dashed, very thick,domain=0.6*180/pi:2.7*180/pi] plot ({e^(0.674394308508077)*cos(\x)},  {e^(0.674394308508077)*sin(\x)}); %

   \draw [black,dashed, very thick,domain=0.6*180/pi:2.7*180/pi] plot ({e^(-0.674394308508077)*cos(\x)},
   {e^(-0.674394308508077)*sin(\x)});%




\node at (0,0) [right] {$0$}; %
\node at (0,0)  {$\bullet$}; %

\draw [black, very thick] plot [smooth] coordinates
{(1/0.7346390419363942,0) (1.45,0.6) (1.635,1.114)};  %

\draw [black, very thick] plot [smooth] coordinates
{(0.7346390419363942,0) (1.65/2.4625,0.53/2.4625) (1.635/3.914221,1.114/3.914221)};  %



\draw [black, very thick] plot [smooth] coordinates
{({cos(1.4*180/pi)}, {sin(1.4*180/pi)}) (0.05,1.24) (-0.1,1.6) (-0.18,1.95)};  %

\draw [black, very thick] plot [smooth] coordinates
{({cos(1.4*180/pi)}, {sin(1.4*180/pi)}) (0.05/1.5401,1.24/1.5401) (-0.1/2.57,1.6/2.57) (-0.18/3.8349,1.95/3.8349)};  %

\node at ({cos(1.4*180/pi)-0.05}, {sin(1.4*180/pi)+0.15}) [right] {$z_0$}; %

\node at (0.72,0.7) [right] {$L_1$}; %

\node at (-0.75,1.04) [right] {$S_1$}; %

\node at (-0.05,0.7) [right] {$L'_1$}; %

\node at (-0.13,1.6) [right] {$L''_1$}; %

\node at (-1.3,0.9) [right] {$G_z$};  %

\node at (-0.67,0.6) [right] {$G_z^-$};
\end{tikzpicture} %
\end{minipage} %

\vspace{0cm} %
\caption{Trajectory structure in the case (b).} 
\end{figure}



Since the quadratic differentials $Q_-(z)\,dz^2$ and
$Q_+(z)\,dz^2$ were extended to $\overline{\mathbb{C}}$ via
reflection principle they are symmetric with respect to the unit
circle $\mathbb{T}$  and can be expressed in terms of rational
functions as follows. Let
$B_-^0(z)=\prod_{k=1}^{n_0^-}(z-c_k)(1-\bar c_k z)$ and
$B_-^\infty(z)=\prod_{k=1}^{n_\infty^-}(z-p_k)(1-\bar p_k z)$,
where the products are taken  over all zeros (counting
multiplicity) of $Q_-(z)\,dz^2$ in the unit disk $\mathbb{D}$ and
over all poles (counting multiplicity) of $Q_-(z)\,dz^2$ in the
unit disk $\mathbb{D}$, respectively.  Also, let
$P_-^0(z)=\prod_{k=1}^{m_0^-}
e^{-i(\pi+\alpha_k)/2}(z-e^{i\alpha_k})$ and
$P_-^\infty(z)=\prod_{k=1}^{m_\infty^-} e^{-i(\pi
+\beta_k)/2}(z-e^{i\beta_k})$, where the products are taken over
all zeros (counting multiplicity) of $Q_-(z)\,dz^2$ on the unit
circle $\mathbb{T}$ and over all poles (counting multiplicity) of
$Q_-(z)\,dz^2$ on the unit circle $\mathbb{T}$, respectively. The
coefficients $e^{-i(\pi+\alpha_k)/2}$ and $e^{-i(\pi+\beta_k)/2}$
in the products representing $P_-^0(z)$ and $P_-^\infty(z)$ are
chosen in this form to have a real coefficient $C_-$ in the
expression for $Q_-(z)\,dz^2$ in the formula (\ref{2.4.4}) below.
Let $B_+^0(z)$, $B_+^\infty(z)$, $P_+^0(z)$, and $P_+^\infty(z)$
denote similar products for the quadratic differential
$Q_+(z)\,dz^2$ and let $n_0^+$, $n_\infty^+$, $m_0^+$,
$m_\infty^+$ denote the number of terms in the corresponding
products.

In the notations introduced above, the quadratic
differentials $Q_-(z)\,dz^2$ and $Q_+(z)\,dz^2$ can be written as %
\begin{equation} \label{2.4.4} %
Q_-(z)\,dz^2=C_-\frac{P_-^0(z)B_-^0(z)}{P_-^\infty(z)B_-^\infty(z)}\,dz^2,
\quad \quad
Q_+(z)\,dz^2=C_+\frac{P_+^0(z)B_+^0(z)}{P_+^\infty(z)B_+^\infty(z)}\,dz^2
\end{equation}
with some  nonzero constants $C_-$ and $C_+$.

In fact, the constants $C_-$ and $C_+$ are real. Indeed,  it
follows from a well known formula linking the number of zeros and
the number of poles of a quadratic differential on
$\overline{\mathbb{C}}$ (see Lemma~3.2 in \cite{J2}), that %
\begin{equation}  \label{2.4.6} %
2n_\infty^- +m_\infty^- -2n_0^--m_0^-=4 \quad {\mbox{and}}  \quad
2n_\infty^+ +m_\infty^+ -2n_0^+-m_0^+=4.
\end{equation} %
Next, to show that $C_-$ is real, we remind that the unit circle
$\mathbb{T}$ consists of arcs of trajectories and/or arcs of
orthogonal trajectories of the quadratic differential
$Q_-(z)\,dz^2$ and therefore $\arg(Q_-(e^{i\theta})\,dz^2)=0$
($\mod \pi$) if $z=e^{i\theta}$ is not a critical point of
$Q_-(z)\,dz^2$. The latter equation (after routing calculation
left to the interested reader) shows that $\arg C_-=0$ ($\mod
\pi$). To perform such calculation one may use  the first equation
in (\ref{2.4.4}) and take into account the first equation in
(\ref{2.4.6}) and the fact that the tangent vector $dz$ to the
unit circle $\mathbb{T}$ at $z=e^{i\theta}$ has the form
$dz=ie^{i\theta}$. Same argument shows that $C_+$ is also real.


Combining equation~(\ref{2.4.4}) with Lemma~\ref{Lemma-1},
we obtain the following. %

\begin{theorem} \label{Theorem-3.1.1}
Let $Q(\zeta)\,d\zeta^2$ be a quadratic differential on
$\overline{\mathbb{C}}$ and let $\Gamma$ be a Jordan curve free of
infinite critical points of $Q(\zeta)\,d\zeta^2$ which consists of
a finite number of arcs of trajectories and/or orthogonal
trajectories of $Q(\zeta)\,d\zeta^2$ and their endpoints. Then the
fingerprint $k: \mathbb{T}\to \mathbb{T}$ of $\Gamma$ is given by
a solution to the functional equation (\ref{2.4.3}) with
$\mathcal{A}(z)$ and $\mathcal{B}(z)$ defined by (\ref{2.4.5}) and
$Q_-(z)\,dz^2$ and $Q_+(z)\,dz^2$ defined by (\ref{2.4.4}). %
\end{theorem} %

Theorem~\ref{Theorem-3.1.1} generalizes the first part of
Theorems~\ref{Theorem-2} and \ref{Theorem-3}. The converse
statement for this theorem, similar to the converse statements of
Theorems~\ref{Theorem-2} and \ref{Theorem-3}, will require
additional restrictions which are discussed below.

Suppose that $Q_-(z)\,dz^2$ and $Q_+(z)\,dz^2$ are quadratic
differentials satisfying conditions (\ref{2.4.4}) and
(\ref{2.4.6}). Then $\mathbb{T}$ consists of a finite number of
maximal open arcs, $\alpha_1^-,\ldots,\alpha_{k_-}^-$, enumerated
in the counterclockwise direction on $\mathbb{T}$, such that for
every $j$ either $Q_-(z)\,dz^2\ge 0$ for all $z\in \alpha_j^-$ or
$Q_-(z)\,dz^2\le 0$ for all $z\in \alpha_j^-$. The maximality here
means that if one or another inequality holds on $\alpha_j^-$ then
it does not hold on any larger arc containing $\alpha_j^-$. As
before, we allow here the case when $k_-=1$ and
$\alpha_1^-=\mathbb{T}$.

Let $z_j^-=e^{i\theta_j^-}$ denote the initial point of
$\alpha_j^-$, $1\le j\le k_-$, and let
$\alpha_{k_-+1}^-=\alpha_1^-$ and $z_{k_-+1}^-=z_1^-$.
The  points $z_j^-$ can be of two types: %

1) If $Q_-(z)\,dz^2$ has the same sign on arcs $\alpha_{j-1}^-$
and $\alpha_j^-$ then $z_j^-$ is a pole of $Q_-(z)\,dz^2$ of even
order.

 2) If $Q_-(z)\,dz^2$ has different signs on arcs
$\alpha_{j-1}^-$ and $\alpha_j^-$ then $z_j^-$ is a critical point
of  $Q_-(z)\,dz^2$ of odd order.

Let $\alpha_1^+,\ldots,\alpha_{k_+}^+$ and
$z_1^+,\ldots,z_{k_+}^+$ denote similar systems of arcs and points
for the quadratic differential $Q_+(z)\,dz^2$.

\begin{definition} %
We will say that quadratic differentials $Q_-(z)\,dz^2$ and
$Q_+(z)\,dz^2$ satisfying conditions (\ref{2.4.4}) and
(\ref{2.4.6}) are coordinated on $\mathbb{T}$ if  their systems of
arcs $\alpha_1^-,\ldots,\alpha_{k_-}^-$ and
$\alpha_1^+,\ldots,\alpha_{k_+}^+$ have the same number of arcs
(i.e. if $k_-=k_+=k\ge 1$) and if they can be enumerated   in the
counterclockwise direction on $\mathbb{T}$
in  such  a way that the following conditions are satisfied:%
\begin{enumerate}  %
\item[(a)] %
$Q_-(z)\,dz^2\ge 0$ on $\alpha_j^-$ if and only if
$Q_+(z)\,dz^2\ge 0$ on $\alpha_j^+$. %
\item[(b)] %
$z_j^-$ is an infinite critical point of $Q_-(z)\,dz^2$ if and
only if $z_j^+$ is  an infinite critical point of $Q_+(z)\,dz^2$.
\item[(c)] 
If the $Q_-$-length
$|\alpha_j^-|_{Q_-}=\int_{\alpha_j^-}\sqrt{Q_-(z)}\,dz$ of
$\alpha_j^-$ is finite then the $Q_+$-length
$|\alpha_j^+|_{Q_+}=\int_{\alpha_j^+}\sqrt{Q_+(z)}\,dz$ of
$\alpha_j^+$ is also finite and
$|\alpha_j^-|_{Q_-}=|\alpha_j^+|_{Q_+}$.
\end{enumerate} %
\end{definition} %

We emphasize here that the arcs $\alpha_j^-$ and $\alpha_j^+$ may
contain zeroes of even order of the quadratic differentials
$Q_-(z)\,dz^2$ and $Q_+(z)\,dz^2$, respectively, and that the
numbers of such  zeros contained in $\alpha_j^-$ and $\alpha_j^+$
may be different.

\begin{theorem} \label{Theorem-4} %
Suppose that the  quadratic differentials $Q_-(z)\,dz^2$ and
$Q_+(z)\,dz^2$ are coordinated on $\mathbb{T}$.  Then there is a
quadratic differential $Q(\zeta)\,d\zeta^2$ defined on
$\overline{\mathbb{C}}$ and a closed Jordan curve  $\Gamma$
consisting of arcs of trajectories and/or orthogonal trajectories
of $Q(\zeta)\,d\zeta^2$  such that the  quadratic differentials
$Q_-(z)\,dz^2$ and $Q_+(z)\,dz^2$ are pullbacks of the quadratic
differential $Q(\zeta)\,d\zeta^2$ under the appropriate mappings
$\varphi_-$ and $\varphi_+$ associated with $\Gamma$.
\end{theorem}

\noindent %
\emph{Proof.} The proof of this theorem is rather standard. We
weld (some authors prefer the term ``glue'') the unit disk
$\mathbb{D}$ (with the quadratic differential $Q_-(z)\,dz^2$
defined on it) and its exterior $\mathbb{D}_+$ (with the quadratic
differential $Q_+(z)\,dz^2$ defined on it) along the unit circle
$\mathbb{T}$ in such a way that $Q_-(z)\,dz^2$ and $Q_+(z)\,dz^2$
will be meromorphic extensions of each other across $\mathbb{T}$.

 First, to each
of the arcs $\alpha_j^-$ and $\alpha_j^+$ we assign a
``representative'' $w_j^-=e^{i\beta_j^-}$ and a ``representative''
$w_j^+=e^{i\beta_j^+}$, respectively. Precisely, we put $w_j^-=
z_j^-$,   $w_j^+= z_j^+$ if $z_j^-$ is a finite critical point and
we put $w_j^-= z_{j+1}^-$,   $w_j^+= z_{j+1}^+$ if $z_j^-$ is an
infinite critical point but $z_{j+1}^-$ is a finite critical
point. If both $z_j^-$ and $z_{j+1}^-$ are infinite critical
points of $Q_-(z)\,dz^2$ then representatives
$w_j^-=e^{i\beta_j^-}$ and $w_j^+=e^{i\beta_j^+}$ can be any
points in the arcs $\alpha_j^-$ and $\alpha_j^+$, respectively.
Also, if $\alpha_1^-=\mathbb{T}$ then representatives
$w_j^-=e^{i\beta_j^-}$ and $w_j^+=e^{i\beta_j^+}$ can be any
points of $\mathbb{T}$.  Now the welding procedure can be
performed as follows.

1) The initial   points $z_j^-$ and $z_j^+$ of the arcs
$\alpha_j^-$ and $\alpha_j^+$ are considered as identical as well
as representatives $w_j^-$ and $w_j^+$ of these arcs (if
different).

2) The points $z_-=e^{i\theta_-}\in \alpha_j^-$ and
$z_+=e^{i\theta_+}\in \alpha_j^+$ are considered to be identical
if either $z_-\in l(\theta_j^-,\beta_j^-)$, $z_+\in l(\theta_j^+,\beta_j^+)$ are such that  %
$$ 
\int_{l(\theta_-,\beta_j^-)} \sqrt{Q_-(z)}\,dz=\int_{l(\theta_+,\beta_j^+)}\sqrt{Q_+(z)}\,dz %
$$ 
or  $z_-\in l(\beta_j^-,\theta_{j+1}^-)$, $z_+\in l(\beta_j^+,\theta_{j+1}^+)$ are such that  %
$$ 
\int_{l(\beta_j^-,\theta_-)} \sqrt{Q_-(z)}\,dz=\int_{l(\beta_j^+,\theta_+)}\sqrt{Q_+(z)}\,dz. %
$$ 


 After this identification of boundary points we obtain
a topological sphere $S$. To turn it into a Riemann sphere we
introduce complex structure on $S$ as follows. We assume that the
disk $\mathbb{D}$ and its exterior $\mathbb{D}_+$ serve as charts
for their points and thus the complex parameter for these points
is given by the identity mapping $z=\tau$.

Suppose now that  $z_0\in S$ is defined by identifying the points
$z_-^0=e^{i\theta_-^0}\in \alpha_j^-$ and
$z_+^0=e^{i\theta_+^0}\in \alpha_j^+$, which are not critical
points of the quadratic differentials $Q_-(z)\,dz^2$ and
$Q_+(z)\,dz^2$, respectively. Then  the complex parameter $\tau$
can be assigned as follows. Let $\varepsilon>0$ be small enough
and let $U_\varepsilon=\{\tau: \, |\tau|<\varepsilon\}$. For
$\tau\in U_\varepsilon$, we define a
function $F:U_\varepsilon\to S$ to be the inverse of the function %
\begin{equation} \label{2.4.8} %
\tau=\Phi(z)=\left\{ \begin{array}{ll} %
\int_{z_-^0}^z\sqrt{Q_-(t)}\,dt & {\mbox{if $z\in
\overline{\mathbb{D}}$ is close enough to $z_-^0$}} \\ %
\int_{z_+^0}^z\sqrt{Q_+(t)}\,dt & {\mbox{if $z\in
\mathbb{D}_+$ is close enough to $z_+^0$.}}  %
\end{array} %
\right.
\end{equation} %
We assume here that the branches of the radicals in (\ref{2.4.8})
are chosen such that $\Im \Phi(z)<0$ if $z\in \mathbb{D}$ and $\Im
\Phi(z)>0$ if $z\in \mathbb{D}_+$.  Furthermore, we stress that
the first integral in (\ref{2.4.8}) is taken along pathes in the
unit disk while the second integral is taken along paths in the
exterior of the unit disk. Under these conditions $z=F(\tau)$
defines a conformal parameter in $U_\varepsilon$. Indeed,
(\ref{2.4.8}) implies that $\Phi(z)$ is analytic in a neighborhood
of the point $z_0$ under consideration, possibly except points of
$S$ corresponding to the points of the unit circle. Furthermore,
the identification rule 2) stated above in this proof implies that %
$$ 
\int_{z_-^0}^{z_-} \sqrt{Q_-(t)}\,dt=\int_{z_+^0}^{z_+}\sqrt{Q_+(t)}\,dt %
$$ 
if $z_-=e^{i\theta_0}$ and $z_+=e^{i\theta_+}$ define the same
point of $S$. The latter implies that $\Phi(z)$ is continuous in a
neighborhood of $z_0$. Since $\Phi(z)$ is analytic in a slit
neighborhood of $z_0$ and continuous in the whole neighborhood it
follows that $\Phi(z)$ is analytic on this neighborhood. Also it
follows from (\ref{2.4.8}) that $\Phi'(z_0)\not=0$ if $z_-^0$ is
not a critical point of $Q_-(z)\,dz^2$ and $z_+^0$ is not a
critical point of $Q_+(z)\,dz^2$. Hence, $\Phi(z)$ is one-to-one
in a neighborhood of $z_0$ in this case.

Next, differentiating both parts of  equation (\ref{2.4.8}) and
squaring the resulting equation, we obtain the following equation
for the quadratic differentials:
$$ 
d\tau^2=\left\{ \begin{array}{ll} %
Q_-(z)\,dz^2 & {\mbox{if $z\in
\mathbb{D}$ }} \\ %
Q_+(z)\,dz^2 & {\mbox{if $z\in
\mathbb{D}_+$.}}  %
\end{array} %
\right.
$$ 
The latter equation shows that, in terms of the parameter
$\tau=\Phi(z)$, the quadratic differentials $Q_-(z)\,dz^2$ and
$Q_+(z)\,dz^2$ define the same quadratic differential $d\tau^2$
and therefore they are extensions of each other.

Let $S'$ denote the surface $S$ punctured at the points which
correspond to critical points situated on the unit circle
$\mathbb{T}$ of the quadratic differentials $Q_-(z)\,dz^2$ and/or
$Q_+(z)\,dz^2$. Since $S'$ is a topological sphere with punctures,
which is supplied with a complex structure, it can be mapped
conformally onto a punctured complex sphere (we denote it by
$\mathbb{C}'$) by a one-to-one function $g:S'\to \mathbb{C}'$. Let
$\psi=g^{-1}$ be the inverse of $g$. Since $\psi$ is conformal and
one-to-one on $\mathbb{C}'$, each puncture is a removable singular
point of $\psi$. Thus, $\psi$ can be extended to a conformal
one-to-one mapping from $\overline{\mathbb{C}}$ to $S$. Let
$Q(\zeta)\,d\zeta^2$ be a quadratic differential on $\mathbb{C}'$
defined by %
 \begin{equation} \label{2.4.9} %
Q(\zeta)\,d\zeta^2=\left\{ \begin{array}{ll} %
Q_-(\psi(\zeta))(\psi'(\zeta))^2\,d\zeta^2 & {\mbox{if
$z=\psi(\zeta)\in
\mathbb{D}\cup \mathbb{T}'$}} \\ %
Q_+(\psi(\zeta))(\psi'(\zeta))^2\,d\zeta^2 & {\mbox{if
$z=\psi(\zeta)\in
\mathbb{D}_+\cup \mathbb{T}'$,}}  %
\end{array} %
\right.
\end{equation} %
where $\mathbb{T}'$ denotes the set of points $e^{i\theta}\in
\mathbb{T}$, which are regular for both $Q_-(z)\,dz^2$ and
$Q_+(z)\,dz^2$.

The quadratic differential (\ref{2.4.9}) can be extended to the
whole complex sphere $\overline{\mathbb{C}}$. Indeed, suppose that
$\zeta_0\not\in \mathbb{C}'$. We may assume that $\zeta_0$ is
finite. Then $\zeta_0$ corresponds via the mapping $\psi(\zeta)$
to a critical point $z_-^0\in \mathbb{T}$ of the quadratic
differential $Q_-(z)\,dz^2$ or/and to a critical point $z_+^0\in
\mathbb{T}$ of the quadratic differential $Q_+(z)\,dz^2$. Since
the trajectory structure of $Q(\zeta)\,d\zeta^2$ near
$\zeta=\zeta_0$ is inherited from the quadratic differentials
$Q_-(z)\,dz^2$ and $Q_+(z)\,dz^2$, it follows that $\zeta_0$ is an
infinite critical point of $Q(\zeta)\,d\zeta^2$ if and only if the
points $z_-^0\in \mathbb{T}$ and $z_+^0\in \mathbb{T}$, which
correspond to $\zeta_0$ under the mapping $\psi(\zeta)$, are
infinite critical points of the quadratic differentials
$Q_-(z)\,dz^2$ and $Q_+(z)\,dz^2$, respectively.

Furthermore, $\zeta_0$ is a finite critical point of
$Q(\zeta)\,d\zeta^2$ if and only if at least one of the
corresponding points $z_-^0\in \mathbb{T}$ and $z_+^0\in
\mathbb{T}$ is a finite critical point of the quadratic
differential $Q_-(z)\,dz^2$ or quadratic differential
$Q_+(z)\,dz^2$, respectively.

Summarizing our previous arguments, we conclude that the image of
the unit circle $\Gamma=\psi^{-1}(\mathbb{T})$ is a closed Jordan
curve consisting of a finite number of arcs of trajectories and/or
arcs of orthogonal trajectories of the quadratic differential
$Q(\zeta)\,d\zeta^2$ and the endpoints of these arcs. The latter
completes our proof of Theorem~\ref{Theorem-4}.~\hfill $\Box$ %

\medskip


The critical points of $Q(\zeta)\,d\zeta^2$ situated off $\Gamma$
and the trajectory structure nearby are inherited from the
corresponding critical points of the quadratic differentials
$Q_-(z)\,dz^2$ and $Q_+(z)\,dz^2$.  To understand the trajectory
structure of $Q(\zeta)\,d\zeta^2$ near its critical points on
$\Gamma$, we have to reverse our procedure described in items
(a)-(g) above. For instance, if $Q_-(z)\,dz^2$ has a simple pole
at $z_-^0\in \mathbb{T}$ and $Q_+(z)\,dz^2$ has a zero of order
$2n-1$, $n\ge 1$, at $z_+^0\in \mathbb{T}$ then $\zeta_0$ is a
zero of order $n-1$ of $Q(\zeta)\,d\zeta^2$ (cf. items (c) and
(d)). In all other possible cases for the points $z_-^0$ and
$z_+^0$, the order of the corresponding critical point $\zeta_0$
of the quadratic differential $Q(\zeta)\,d\zeta^2$ can be found in
a similar way by using items (a)-(g).

\medskip


Because  we have a choice of representatives, a quadratic
differential $Q(\zeta)\,d\zeta^2$ of Theorem~\ref{Theorem-4} is
not unique in general. Thus any uniqueness result related to
Theorem~\ref{Theorem-4} will need some additional assumptions as
in our theorem below.

\begin{theorem} \label{Theorem-5} %
Suppose that the assumptions of Theorem~\ref{Theorem-4} are
satisfied. Then, for any given sets of representatives $w_j^-$ and
$w_j^+$, a quadratic differential $Q(\zeta)\,d\zeta^2$ and a
Jordan curve $\Gamma$ as in Theorem~\ref{Theorem-4}  are defined
uniquely up to
a M\"{o}bius transformation. 
\end{theorem} %

\noindent %
\emph{Proof.} Suppose that we have two quadratic differentials
$Q_1(\zeta)\,d\zeta^2$ and $Q_2(\zeta)\,d\zeta^2$ and suppose that
$\Gamma_1$ and $\Gamma_2$ denote the corresponding Jordan curves.
Suppose further that $\varphi_-^1:\mathbb{D}\to \Omega_-^1$,
$\varphi_+^1:\mathbb{D}_+\to \Omega_+^1$,
$\varphi_-^2:\mathbb{D}\to \Omega_-^2$, and
$\varphi_+^2:\mathbb{D}_+\to \Omega_+^2$ denote canonical mappings
onto domains $\Omega_-^1$, $\Omega_+^1$, $\Omega_-^2$, and
$\Omega_-^2$, which  correspond to curves $\Gamma_1$ and
$\Gamma_2$, respectively. We assume here that $Q_-(z)\,dz^2$ is a
pullback of $Q_j(\zeta)\,d\zeta^2$, $j=1,2$, under the mapping
$\varphi_-^j(z)$ and that  $Q_+(z)\,dz^2$ is a pullback of
$Q_j(\zeta)\,d\zeta^2$, $j=1,2$, under the mapping
$\varphi_+^j(z)$. Then the composed mapping %
\begin{equation} \label{2.4.10} %
\varphi(\zeta)=\left\{ \begin{array}{ll} %
\varphi_-^2\left((\varphi_-^1)^{-1}(\zeta)\right) & {\mbox{if
$\zeta\in
\Omega_-^1\cup \Gamma_1$}} \\ %
 \varphi_+^2\left((\varphi_+^1)^{-1}(\zeta)\right) & {\mbox{if $\zeta\in
\Omega_+^1$}}  %
\end{array} %
\right.
\end{equation} %
is conformal on $\Omega_-^1$ and on $\Omega_+^1$. Furthermore, the
quadratic differentials $Q_-(z)\,dz^2$ and $Q_+(z)\,dz^2$ are
pullbacks (under appropriate mappings) of each of the quadratic
differentials $Q_1(\zeta)\,d\zeta^2$ and $Q_2(\zeta)\,d\zeta^2$.
The latter implies that $\varphi(\zeta)$ is continuous on
$\Gamma_1$. Since $\Gamma_1$ is piecewise analytic and
$\varphi(\zeta)$ is continuous on $\Gamma_1$ it follows that
$\varphi(\zeta)$ is conformal on $\Gamma_1$. Thus,
$\varphi(\zeta)$ is conformal and one-to-one on
$\overline{\mathbb{C}}$ and therefore $\varphi(\zeta)$ is a
M\"{o}bius mapping, as required. \hfill $\Box$

\begin{corollary}
Suppose that  the assumptions of Theorem~\ref{Theorem-4} are
satisfied and suppose additionally that the number of arcs
$\alpha_j^-$ is at least two and at least one of the endpoints of
each of the arcs $\alpha_j^-$ is a finite critical point. Then  a
quadratic differential $Q(\zeta)\,d\zeta^2$ and a Jordan curve
$\Gamma$ as in Theorem~\ref{Theorem-4}  are defined uniquely up to
a M\"{o}bius transformation.
\end{corollary} %

\medskip

If both endpoints of $\alpha_j^-$ are infinite critical points of
$Q_-(z)\,dz^2$ and if both $Q_-(z)\,dz^2$ and $Q_+(z)\,dz^2$ have
zeroes on the arcs $\alpha_j^-$ and $\alpha_j^+$ then by choosing
the representatives defined in the proof of
Theorem~\ref{Theorem-4} appropriately we can obtain quadratic
differentials $Q_1(\zeta)\,d\zeta^2$ and $Q_2(\zeta)\,d\zeta^2$,
which have different numbers of distinct zeroes. Thus,
$Q_1(\zeta)\,d\zeta^2$ and $Q_2(\zeta)\,d\zeta^2$ are essentially
different and therefore the generated Jordan curves $\Gamma_1$ and
$\Gamma_2$ are also essentially different. Even in the case when
$\alpha_1^-=\alpha_1^+=\mathbb{T}$ and $\alpha_1^-$ and
$\alpha_1^+$ do not contain critical points we may have
essentially different quadratic differentials for different
choices of representatives. Here are some examples.

\noindent %
\textbf{(1)} In the simplest case when
$$ 
Q_-(z)\,dz^2=Q_+(z)\,dz^2=-\frac{dz^2}{(z+1)^4}, %
$$ 
for any choice of representatives $w_1^-$ and $w_1^+$, the
resulting quadratic differential is
$Q(\zeta)\,d\zeta^2=-\frac{d\zeta^2}{(\zeta+1)^4}$ and $\Gamma$ is
the unit circle (up to a M\"{o}bius transformation). Indeed,
changing variables via $z=\frac{i-\tau}{i+\tau}$, we transform
$Q_-(z)\,dz^2$ into the quadratic differential
$\frac{1}{4}\,d\tau^2$ whose trajectories are horizontal lines.
Then, the procedure of choosing representatives $w_1^-$ and
$w_1^+$ is equivalent to the geometric procedure of sliding the
upper-half plane of a variable $\tau$ with respect to its
lower-half plane. Since these slidings don't affect horizontal
lines and Euclidean lengthes they don't change the trajectory
structure and $Q$-lengths, and therefore the quadratic
differential $\frac{1}{4}\,d\tau^2$ remains unchanged.

The same geometric argument can be applied to identify the form of
the quadratic differential $Q(\zeta)\,d\zeta^2$ in our examples
\textbf{(2)}, \textbf{(3)}, and \textbf{(4)} below. Necessary
technical details in these examples are left to the interested
reader.

\smallskip

\noindent %
\textbf{(2)} In the case when
$$ 
Q_-(z)\,dz^2=Q_+(z)\,dz^2=\frac{dz^2}{(z-1)^2(z+1)^2}, %
$$ 
different choices of pairs of representatives $w_1^-$, $w_2^-$ and
$w_1^+$, $w_2^+$ generate  quadratic differentials
$Q(\zeta)\,d\zeta^2=e^{i\alpha}\frac{d\zeta^2}{(\zeta-1)^2(\zeta+1)^2}$
with $-\pi<\alpha<\pi$ depending on the choice of representatives.
In case $\alpha=0$, the resulting curve $\Gamma$ coincides with
the unit circle and in the case $\alpha\not=0$ it  consists of a
pair of logarithmic spirals with foci at $\pm 1$; see Figure~3.

\smallskip


\begin{figure}\label{fig-2}

$$\includegraphics[scale=.5,angle=0]{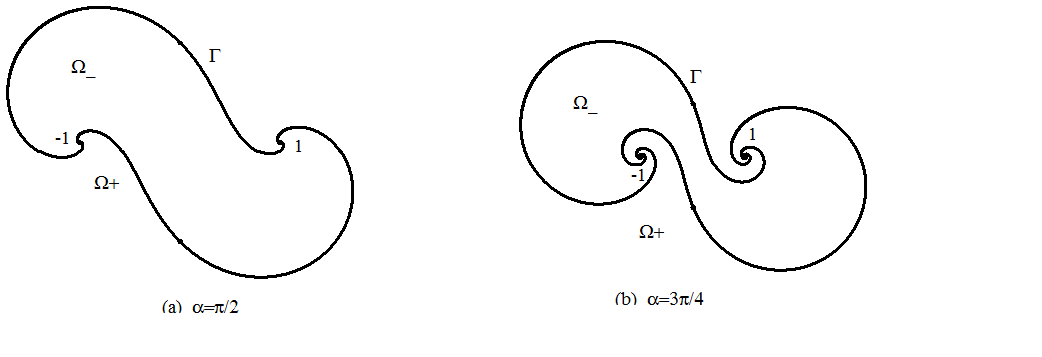}
$$
\vspace{-1cm} %
\caption{$\Gamma$ consisting of two spirals with different $\alpha$.} %
\end{figure}

\FloatBarrier

\noindent %
\textbf{(3)} Let %
$$ 
Q_-(z)\,dz^2=Q_+(z)\,dz^2=\frac{(z-1)^2}{(z+1)^6}\,dz^2.
$$ 
For this pair of quadratic differentials we have only two
possibilities. The first one occurs when  $w_1^-=w_1^+=1$. In this
case
$Q(\zeta)\,d\zeta^2=\frac{(\zeta-1)^2}{(\zeta+1)^6}\,d\zeta^2$ and
$\Gamma$ is the unit circle. For any other choice of
representatives $w_1^-$ and $w_1^+$, we have quadratic
differentials of the form $Q(\zeta)\,d\zeta^2=C
\frac{\zeta^2+1}{(\zeta+1)^6}\,d\zeta^2$ with some nonzero $C\in
\mathbb{C}$ depending on the choice of representatives. In this
case $\Gamma$ consists of three critical trajectories as it is
shown in Figure~4.


\begin{figure} \label{fig-3}

$$\includegraphics[scale=.5,angle=0]{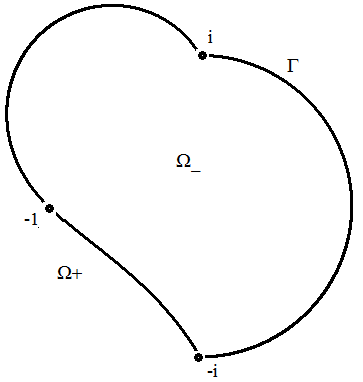}
$$
\vspace{-1cm} %
\caption{$\Gamma$ consisting of three critical trajectories.} %
\end{figure}


\smallskip

\noindent %
\textbf{(4)} For a given $p$, $0<p<1$, consider quadratic differentials %
$$ 
Q_-(z)\,dz^2=Q_+(z)\,dz^2=\frac{(z-p)(z-1/p)}{z^3}\,dz^2.
$$ 
In this case, for every $p$, we have a family of quadratic
differentials of the form %
$$ %
Q(\zeta)\,d\zeta^2 =e^{i\alpha}
\frac{(\zeta-\rho)(\zeta-e^{i\beta}/\rho)}{\zeta^3}\,d\zeta^2
$$ %
 with $0<\rho<1$ and  $0\le \beta<\pi$ depending on $p$ and on the choice of
representatives and with
$$ %
\alpha=-\arg\left(i\int_0^{2\pi}\sqrt{(e^{i\theta}-\rho)(e^{i\theta}-e^{i\beta}/\rho)}\,
e^{i\theta/2}\,d\theta \right). %
$$ %

For $\rho=1/3$ and for three different choices of $\beta$, the
critical trajectories of corresponding quadratic differentials and
the resulting curves $\Gamma$ are shown in  Figure~5.

We note here that even when quadratic differentials have a small
number of critical points they still can generate a rich variety
of topological structures depending on the positions of these
critical points; see, for example, Section~11 in \cite{ShSol},
which presents all possible structures of critical trajectories
for a quadratic differential
$Q(z)\,dz^2=e^{i\alpha}\frac{(z-p_1)(z-p_2)}{(z^2-1)^2}\,dz^2$.

Two more  examples will be given in Section~5, where we discuss
polygonal curves. In the latter case, the corresponding quadratic
differentials $Q_-(z)\,dz^2$ and $Q_+(z)\,dz^2$ are related to the
well-known  Schwarz-Christoffel integrals.

\FloatBarrier

\begin{figure} \label{fig-4}

$$\includegraphics[scale=.4,angle=0]{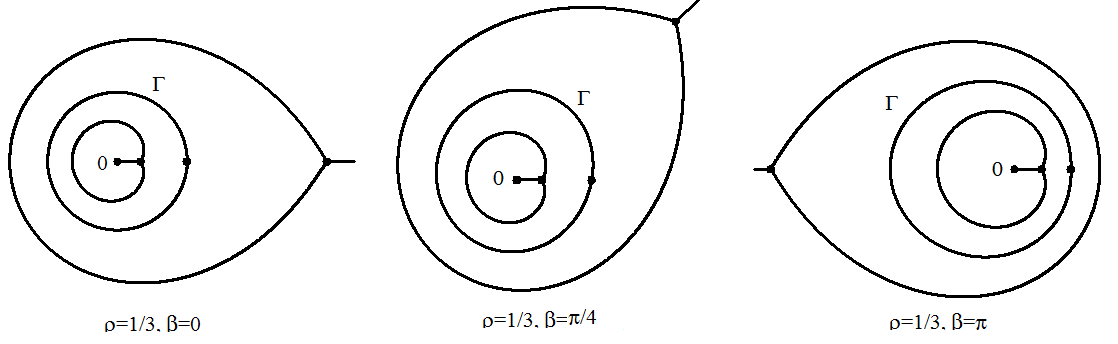}
$$
\vspace{-0.8cm} %
\caption{$\Gamma$ consisting of one regular trajectory.} %
\end{figure}


\section{Modulated graphs and quadratic differentials}

\setcounter{equation}{0} %

In this section we discuss a procedure designed for construction
of quadratic differentials with required properties. We note that
the full generality of results presented here is not required for
the rest of this paper but they can be useful for the future work
on image recognition. For the purposes of this section, it is
convenient to use the terminology of  graph theory, although we do
not need any deep results from it. Our exposition here follows the
lines of
 the paper \cite{S2}  with the necessary changes to
include doubly-connected domains. Let $\mathrm{G}^0=(V,E,F)$ be a
graph embedded in $\overline{\mathbb{C}}$ with the set of vertices
$V=\{v_k\}$, the set of edges $E=\{e_{kj}^i\}$, and the set of
faces $F=\{f_k\}$. We will assume that edges are enumerated in
such a way that $e_{kj}^i$ separates faces $f_k$ and $f_j$. Then,
of course, $e_{kj}^i=e_{jk}^i$. In case $k=j$, an edge $e_{kk}^i$
separates $f_k$ from itself. In this case $e_{kk}^i$ is called a
\emph{cut edge}.

We suppose further that $\mathrm{G}^0$ is \emph{at most doubly
connected}. The latter means that each face is either simply
connected or a doubly connected domain. For convenience, we assume
additionally that at most doubly connected graphs do not have
connected components which are cyclic graphs.

If $f_k$ is simply connected then its boundary $\partial f_k$
represents a closed walk in $\mathrm{G}^0$ in which each vertex
may be visited more than once and each cut edge is traversed
twice. To be definite, we assume that all boundary walks
considered in this paper are oriented in a counterclockwise
direction with respect to corresponding faces $f_k$. If $f_k$ is
doubly connected then its boundary consists of two boundary walks.
To fix notation for future use, one of these boundary walks will
be denoted by $\partial_1 f_k$ and the other one by $\partial_2
f_k$. We note here that, by our definition of at most doubly
connected graphs, all boundary walks are distinct. Then every walk
along each boundary component of $f_k$ contains at least one
vertex which order is not two.



Next, we will load $\mathrm{G}^0$ with three sets, the set of
weights $W=\{w_{kj}^i\}$, the set of heights $H=\{h_k\}$, and the
set of twists $T=\{t_k\}$.
The element $w_{kj}^i$ of $W$ assigns positive weight to the edge
$e_{kj}^i$. Weighted graphs are standard in the graph theory.
Using weights, we can define lengthes  (or total weights)
$\alpha_k$ of boundary
walks around simply connected faces $f_k$ as follows: %
\begin{equation} \label{3.1} %
\alpha_k=\sum_{i,j} w_{kj}^i,
\end{equation} %
where the weight $w_{kk}^i$ of each cut edge $e_{kk}^i\subset
\partial f_k$ is counted twice. If $f_k$ is doubly connected then
lengthes $\alpha_k^1$ of $\partial_1 f_k$ and $\alpha_k^2$ of
$\partial_2 f_k$ are defined in a similar way: %
\begin{equation} \label{3.2} %
\alpha_k^m=\sum w_{kj}^i, \quad m=1,2, %
\end{equation} %
where the sum is taken over all weights $w_{kj}^i$ corresponding
to the edges $e_{kj}^i\subset \partial_m f_k$, $m=1,2$. As before
the weight $w_{kk}^i$ of each cut edge $e_{kk}^i\subset
\partial_m f_k$ is counted twice.
An at most doubly connected graph is called \emph{coordinated} if
$\alpha_k^1=\alpha_k^2$ for each  doubly connected face~$f_k$. In
the case of coordinated graphs we will write $\alpha_k$ instead of
$\alpha_k^1$ and $\alpha_k^2$.

Two other assigned sets, $H=\{h_k\}$ with $h_k>0$  and
$T=\{t_k\}$, are not standard in the graph theory. These sets are
assigned to doubly connected faces $f_k$ only with the purpose to
introduce a sort of a \emph{moduli space structure} on graphs.
Specifically, the height $h_k>0$ can be thought as a
(combinatorial) distance between boundary components $\partial_1
f_k$ and $\partial_2 f_k$. Then the quotient %
\begin{equation}
m(f_k)=\frac{h_k}{\alpha_k}
\end{equation} %
can be considered as a (combinatorial) module of a doubly
connected face $f_k$ (compare this with the formulas (\ref{3.3})
and (\ref{3.3.1}) below).


To explain 
what the term ``twist of $f_k$'' means, we   represent the
boundary walks $\partial_1 f_k$ and $\partial_2 f_k$ by cyclic
graphs $\gamma_{k,1}$ and $\gamma_{k,2}$ with vertices and edges
on the unit circle. Any such representation has the form %

\begin{equation} \label{3.02} %
\gamma_{k,m}=c_1(k,m) l_1(k,m) 
\ldots c_N(k,m)l_N(k,m) c_{N+1}(k,m), \quad m=1,2,
\end{equation}  %
where  
\begin{equation} \label{3.05} %
c_s(k,m)= e^{i\theta_s(k,m)}, \quad \quad
l_s(k,m)=\{e^{i\theta}:\, \theta_s(k,m)<\theta<\theta_{s+1}(k,m)\}
\end{equation}  %
and $c_{N+1}(k,m)=e^{i\theta_{N+1}}=e^{i(\theta_1+2\pi)}$.  In
these formulas, $N=N(k,m)$ denotes the number of edges in the
corresponding walk, where each cut edge is counted twice. In a
walk along a boundary component of $f_k$ some vertices can be
visited several times. Hence, the same vertex of the graph $G^0$
can be represented by more than one of the vertices $c_s(k,m)$ of
the string (\ref{3.02}). Furthermore, we will  assume that the
lengths of the arcs $l_s(k,m)$ are proportional to the weights
$w_{kj}^i$ of the edges $e_{kj}^i$, which are
represented by these arcs. Thus, we assume that %
\begin{equation} \label{3.07} %
\frac{1}{2\pi}
\left(\theta_{s+1}(k,m)-\theta_s(k,m)\right)=\frac{w_{kj}^i}{\alpha_k}
\end{equation}
if the arc $l_s(k,m)$ represents the edge $e_{kj}^i\in
\partial_m f_k$.

A pair $(\gamma_{k,1},\gamma_{k,2})$ of cyclic graphs defined by
(\ref{3.02}) -- (\ref{3.07}) will be called a representation of
the pair of boundary walks $(\partial_1 f_k,\partial_2 f_k)$. Two
such representations $(\gamma'_{k,1},\gamma'_{k,2})$ and
$(\gamma''_{k,1},\gamma''_{k,2})$ are equivalent if  one of them
is a rotation of the other around the origin; i.e. if there is
real $\alpha$ such that $\gamma''_{k,1}=e^{i\alpha}\gamma'_{k,1}$
and $\gamma''_{k,2}=e^{i\alpha}\gamma'_{k,2}$.

\begin{definition} %
By a twist $t_k$ of $f_k$ we mean the class of equivalent
representations  $(\gamma_{k,1},\gamma_{k,2})$ of the pair of
boundary walks $(\partial_1 f_k,\partial_2 f_k)$.
\end{definition}

Thus, twist $t_k$ characterizes the relative positions of cyclic
graphs $\gamma_{k,1}$ and $\gamma_{k,2}$ and is not a pure
numerical characteristic. To assign a numerical value to $t_k$, we
may fix a position of one vertex for each of the cyclic graphs
$\gamma_{k,1}$ and $\gamma_{k,2}$.  For instance, we may choose %
\begin{equation} \label{3.06} %
c_1(k,m)=1 \quad c_1(k,2)=e^{i\tau_k} \quad {\mbox{with $0\le
\tau_k <2\pi$.}} %
\end{equation} %

Then, the parameter $\tau_k$  is  the angle which characterizes
rotation of the cyclic graph $\gamma_{k,2}$ with respect to the
cyclic graph $\gamma_{k,1}$ but the value of this angle  depends
on the choice of a pair of ``initial'' vertices on the graphs
$\gamma_{k,1}$ and $\gamma_{k,2}$.

\begin{definition} %
A coordinated at most doubly connected   graph
$\mathrm{G}^0=(V,E,F)$ loaded with the set of weights
$W=\{w_{kj}^i\}$ assigned to its edges $e_{kj}^i$ and  sets of
heights $H=\{h_k\}$ and  twists $T=\{t_k\}$ assigned to its doubly
connected faces $f_k$ will be called a
modulated graph and will be denoted by $\mathrm{G}=(V,E,F,W,H,T)$. %
\end{definition} %

\bigskip

Next, we will discuss graphs associated with quadratic
differentials. Let $Q(z)\,dz^2$ be a quadratic differential on
$\overline{\mathbb{C}}$.  Consider sets $V_Q$, $E_Q$, and $F_Q$,
where $V_Q$ is the set of all critical points of $Q(z)\,dz^2$,
$E_Q$ is  the set of all critical trajectories of $Q(z)\,dz^2$,
and $F_Q$ is the set consisting of all open components of
$\overline{\mathbb{C}}\setminus \overline{E_Q}$. It follows from
Jenkins' \emph{Basic Structure Theorem} (see Theorem~3.5 in
\cite{J2}) that $F_Q$ consists of a finite number of \emph{circle
domains}, \emph{ring domains}, \emph{strip domains}, and \emph{end
domains} of the quadratic differential $Q(z)\,dz^2$. A \emph{ring
domain} in this context means a doubly connected domain.
Furthermore, the interior of $\overline{E_Q}$ may have a finite
number of connected component, which are called \emph{density
domains}.

For the purposes of this paper, we need  quadratic differentials
with circle domains and ring domains only; i.e. without strip
domains, end domains,  and density domains. Quadratic
differentials with these properties are commonly known as
\emph{Jenkins-Strebel} quadratic differentials.

Suppose now that $Q(z)\,dz^2$ is a Jenkins-Strebel quadratic
differential. Then the triple  $(V_Q,E_Q,F_Q)$ can be considered
as a graph, which is at most doubly connected in our terminology.
Furthermore, the $Q$-length of every critical trajectory of a
Jenkins-Strebel quadratic differential is finite. Hence, we may
suppose further that every critical trajectory $e_{kj}^i\in E_Q$
carries a weight $w_{kj}^i$ equal to the $Q$-length of this
critical trajectory. The set of all these weights will be denoted
by $W_Q$. The string of four set $V_Q$, $E_Q$, $F_Q$, and $W_Q$
defines the \emph{weighted critical graph}
$\mathrm{G}_Q^0=(V_Q,E_Q,F_Q,W_Q)$ of the quadratic differential
$Q(z)\,dz^2$.

We note here that the weighted critical graph of a Jenkins-Strebel
quadratic differential is at most doubly connected and
coordinated.

 Information provided by the weighted graph
$\mathrm{G}_Q^0$ is usually sufficient  when we work with
Jenkins-Strebel quadratic differentials which faces are simply
connected domains. To have graphs carrying more information about
Jenkins-Strebel quadratic differentials with doubly connected
faces we will load their weighted graphs with two additional sets,
the set of heights $H_Q$ and the set of twists $T_Q$. As we
mentioned above, heights and twists are defined  for ring domains
of $Q(z)\,dz^2$ only.

Let $f_k\in F_Q$ be a ring domain of $Q(z)\,dz^2$. By the height
$h_k$ of $f_k$ we mean the $Q$-distance between boundary
components of $f_k$. The latter coincides with the $Q$-length of
an arc of any orthogonal trajectory, which joins the boundary
components of $f_k$. Thus, $0<h_k<\infty$. Now we can define the
set of heights of $Q(z)\,dz^2$ as $H_Q=\{h_k\}$.

Let $m(f_k)$ denote the module of a ring face $f_k\in F_Q$ with
respect to the family of
curves separating boundary components of $f_k$. Then %
\begin{equation} \label{3.3} %
m(f_k)=\frac{h_k}{\alpha_k},
\end{equation} %
where $\alpha_k$ denotes the length of the boundary walk around a
boundary component of the face $f_k$. Formula~(\ref{3.3}) relates
the set of \emph{weights} and the set of \emph{heights} with the
set of \emph{moduli} of doubly connected faces $f_k\in F_Q$. The
latter explains our choice of the term ``modulated graphs''.

 To introduce the twist $t_k\in T_Q$ of a ring face $f_k\in F_Q$, we
will consider  a canonical conformal mapping $\Phi_k(z)$ from
$f_k$ onto a circular annulus $A(1,R_k)$, where $A(1,R)=\{z:\,
1<|z|<R\}$.
The latter mapping is  given by %
$$ 
\Phi_k(z)=\exp\left(-\frac{2\pi i}{\alpha_k} \int_{z_0}^z
\sqrt{Q(\tau)}\,d\tau\right), \quad \quad z\in f_k.
$$ 
We assume here that the initial point of integration, $z_0$,
belongs to $\partial_1f_k$ and that the branch of the radical is
chosen such that $|\Phi_k(v_{k_2})|>1$. Since the module of a
doubly connected domain is conformally invariant and
$m(A(1,R))=\frac{1}{2\pi}\log R$ we conclude from (\ref{3.3}) that %
\begin{equation} \label{3.3.1} %
\frac{h_k}{\alpha_k}=\frac{1}{2\pi}\log R_k.
\end{equation} %

The function $\Phi_k(z)$ maps the boundary components
$\partial_1f_k$ and $\partial_2 f_k$ onto circles $C_1$ and
$C_{R_k}$. Moreover, this mapping is continuous and one-to-one in
the sense of boundary correspondence. Therefore, the normalized
function $\varphi_k(z)=\frac{\Phi_k(z)}{|\Phi_k(z)|}$ maps
boundary walks $\partial_1 f_k$ and $\partial_2 f_k$ onto cyclic
graphs, call them $\gamma_{k,1}^Q$ and $\gamma_{k,2}^Q$, placed on
the unit circle. Now, the pair $(\gamma_{k,1}^Q,\gamma_{k,2}^Q)$
of these cyclic graphs defines a twist $t_k^Q$ which will be
called the twist of the ring domain $f_k$ (of the quadratic
differential $Q(z)\,dz^2$). Let $T_Q=\{t_k\}$ denote the set of
twists of the quadratic differential $Q(z)\,dz^2$.

\begin{definition} %
The weighted  graph $\mathrm{G}_Q^0=(V_Q,E_Q,F_Q,W_Q)$  of the
quadratic differential $Q(z)\,dz^2$ loaded with two additional
sets, the set of heights $H_Q=\{h_k\}$ and the set of twists
$T_Q=\{t_k\}$ assigned as above to  ring domains of $Q(z)\,dz^2$,
will be called the
critical modulated graph of $Q(z)\,dz^2$ and will be denoted by $\mathrm{G}_Q=(V_Q,E_Q,F_Q,W_Q,H_Q,T_Q)$. %
\end{definition} %


Before we state the main result of this section we recall the
following. A function $\tau:\overline{\mathbb{C}}\times[0,1]\to
\overline{\mathbb{C}}$ is called a deformation of
$\overline{\mathbb{C}}$  if the following conditions are
fulfilled: %

1)  $\tau$  is continuous on $\overline{\mathbb{C}}\times [0, 1]$, %

2) $\tau(\cdot,0)$ is the identity mapping, %

3) for every fixed $t\in [0, 1]$, $\tau(\cdot,t)$  is a
homeomorphism from $\overline{\mathbb{C}}$  onto itself.

\medskip

Every deformation $\tau$ transforms a given modulated graph $G$
embedded in $\overline{\mathbb{C}}$ into a \emph{homeomorphic}
modulated graph $G_\tau$ embedded in $\overline{\mathbb{C}}$. We
assume additionally that deformation respect the weights of edges
of $G$ and the heights and twists of doubly connected faces of
$G$.

\begin{theorem} \label{Theorem 3.1}
For every at most doubly connected modulated graph %
$$ %
\mathrm{G}=\{V,E,F,W,H,T\} %
$$ %
 embedded into $\overline{\mathbb{C}}$
there exists a Jenkins-Strebel quadratic differential $Q(z)\,dz^2$
on $\overline{\mathbb{C}}$, the critical modulated graph of which
is homeomorphic to $\mathrm{G}$ on $\overline{\mathbb{C}}$.

Furthermore, $Q(z)\,dz^2$ is defined uniquely  up to a M\"{o}bius
automorphism of $\overline{\mathbb{C}}$.
\end{theorem}

This theorem is also valid  for quadratic differentials and
modulated graphs embedded into an arbitrary compact Riemann
surface; see \cite{S2}, \cite{Em}. The study of weighted graphs
and quadratic differentials on compact Riemann surfaces with
relation to extremal problems was initiated by this author and for
graphs with simply connected faces the result stated above was
proved in Theorem~1 in \cite{S2}. Then an extension of this
theorem to include quadratic differentials with doubly connected
domains with essentially the same proof was given by E.~Emel'yanov
who used a slightly different terminology; see \cite{Em}.

\section{Lemniscates and  quadratic differentials}

\setcounter{equation}{0}

Let $f(z)$ be a nonconstant meromorphic function on a domain
$D\subset\overline{ \mathbb{C}}$. For $0\le
c\le \infty$, the lemniscate of $f(z)$ at level $c$ is defined by equation %
$$ 
L_f(c)=\{z\in D: \, |f(z)|=c\}.
$$ 

Every connected component of $L_f(c)$ will be called a
\emph{lemniscate component}. Thus,   $L_f(c)$ is just a level set
of $f(z)$ with value $c$. Since $L_f(0)$ is the set of zeroes of
$f(z)$ and $L_f(\infty)$ is the set of its poles, each lemniscate
component of $L_f(0)$ and $L_f(\infty)$ is a point.

Suppose now that $0<c<\infty$ and that $L_f(c)$ is not empty. Then
each lemniscate component of $L_f(c)$ is either a closed Jordan
analytic curve or it consists of Jordan analytic arcs and their
endpoints in $D$.

Let $L'_f(c)$ be a lemniscate component having a tangent vector
$dz$ at its point $z$. The gradient of the real valued function
$\log|f(z)|$ at $z$ can be calculated as follows: %
$$ 
{\mbox{grad}}( \log|f(z)|)= 2\frac{\partial}{\partial
\overline{z}}\log|f(z)|=\overline{\left(\frac{f'(z)}{f(z)}\right)}.
$$ 
Hence, the tangent vector $dz$ to $L'_f(c)$ at $z$ is given by %
$$ 
dz=i \overline{\left(\frac{f'(z)}{f(z)}\right)}.
$$ 
Multiplying both sides of this equation by
$\frac{i}{2\pi}\frac{f'(z)}{f(z)}$ and then squaring, we obtain %
\begin{equation} \label{4.4} %
-\frac{1}{4\pi^2}\left(\frac{f'(z)}{f(z)}\right)^2\,dz^2=\frac{1}{4\pi^2}\,
\left|\frac{f'(z)}{f(z)}\right|^4.
\end{equation} %
The left hand-side of (\ref{4.4}) defines a quadratic
differential, which will be denoted by $Q_f(z)\,dz^2$; i.e., %
\begin{equation} \label{4.5} %
Q_f(z)\,dz^2=-\frac{1}{4\pi^2}\left(\frac{f'(z)}{f(z)}\right)^2\,dz^2
\end{equation} %
Now equation (\ref{4.4}) shows that %
$$ 
Q_f(z)\,dz^2>0, %
$$ 
if $dz$ is a tangent vector to the lemniscate of $f(z)$ passing
through $z$.

From here we conclude that, for $0<c<\infty$, each lemniscate
component is a trajectory of the quadratic differential
(\ref{4.5}). Furthermore, each critical point of $Q_f(z)\,dz^2$
has even order and that endpoints of the analytic arcs
constituting lemniscate components are zeroes of $f'(z)$, which
are not zeroes of $f(z)$.

Since the quadratic differential $Q_f(z)\,dz^2$ is generated by
the logarithmic derivative its trajectory structure is  not of a
general form.  Let $z_0\not=\infty$ be a zero or pole of $f(z)$ of
order $s\ge 1$. Then the logarithmic derivative  $f'(z)/f(z)$ can
be represented
as  %
\begin{equation} \label{4.7} %
\frac{f'(z)}{f(z)}=\pm \frac{s}{z-z_0}+{\mbox{higher powers of
$(z-z_0)$}}
\end{equation} %
with ``$+$'' sign for zeros and ``$-$'' sign for poles. Similarly,
if $z_0=\infty$ then %
\begin{equation} \label{4.7.1} %
\frac{f'(z)}{f(z)}=\pm \frac{s}{z}+{\mbox{higher powers of $z$}}
\end{equation} %
with ``$-$'' sing for zeros and ``$+$'' sign for poles.
 In what follows we will assume that
$z_0\not=\infty$, the case $z_0=\infty$ requires only minor
notational  changes.

It follows from (\ref{4.7}) that the quadratic differential
$Q_f(z)\,dz^2$ has the following representation near $z_0$: %
$$ 
Q_f(z)\,dz^2=-\frac{s^2}{4\pi^2}\frac{1}{(z-z_0)^2}+{\mbox{higher
powers of $(z-z_0)$.}} %
$$ 
The latter implies that every zero or pole of $f(z)$ represents a
second order pole of the quadratic differential $Q_f(z)\,dz^2$
with circular structure of trajectories near $z_0$. Moreover, the
$Q_f$-length of each trajectory $\gamma$ separating $z_0$ from
the boundary of $D$ and other critical points is %
$$ 
\left|\gamma\right|_{Q_f}=\int_\gamma |Q_f(z)|^{1/2}\,|dz|=s. %
$$ 

Our discussion above implies, in particular, that the trajectory
structure of $Q_f(z)\,dz^2$ does not include end domains and strip
domains having poles on their boundaries in $D$. Also, the
trajectory structure of $Q_f(z)\,dz^2$ does not include density
domains since otherwise some level set of $f(z)$ would be dense in
such a domain. Then, $|f(z)|$ must be constant on a density domain
and therefore $f(z)$ must be constant on $D$ contradicting our
assumption.

Suppose now that $\gamma$ is a closed trajectory of $Q_f(z)\,dz^2$
and that $f(z)$ is meromorphic on a connected component $D_0$ of
$\overline{\mathbb{C}}\setminus \gamma$.  Then the $Q_f$-length of
$\gamma$ can be calculated as follows: %

\begin{equation} \label{4.8.1} %
|\gamma|_{Q_f}=\int_\gamma
\sqrt{Q_f(z)}\,dz=\frac{1}{2\pi}\left|\int_\gamma
\frac{f'(z)}{f(z)}\,dz\right| =|N_P-N_Z|,
\end{equation} %
where $N_P$ and $N_Z$ denote the number of poles and number of
zeros of $f(z)$ in $D_0$ (counting multiplicity), respectively.

Formula (\ref{4.8.1}) shows, in particular, that a regular
trajectory of $Q_f(z)\,dz^2$ cannot enclose an equal number of
poles and zeros of $f(z)$ (if $f(z)$ is meromorphic on $D_0$). In
terms of lemniscates, this fact can be stated as follows. Suppose
that a Jordan curve $\gamma$ is a subset of some lemniscate
$L_f(c)$. Suppose further that $f(z)$ is meromorphic on a
component $D_0$ of $\overline{\mathbb{C}}\setminus \gamma$ and
$N_P=N_Z$. Then $\gamma$ contains a critical point of $f(z)$ or,
in other words, $Q_f(z)\,dz^2$ has a zero on $\gamma$.

Our previous discussions in this section deal with meromorphic
functions on an arbitrary domain on the complex sphere. All these
results and formulas can be extended to the case of meromorphic
functions defined on Riemann surfaces.

In the rest of this section, we will focus specifically on the
case when $D=\overline{\mathbb{C}}$.  In this case, every
meromorphic function is a rational function and therefore we will
denote such a function by $R(z)$ rather than $f(z)$. Let
$\Gamma\subset \mathbb{C}$ be a connected component of $L_R(c)$
which does not contain critical points of $R(z)$. Let $N_P^-$ and
$N_Z^-$ denote the number of poles and zeros of $R(z)$ in the
domain $\Omega_-$ and let $N_P^+$ and $N_Z^+$ denote the number of
poles and zeros of $R(z)$ in the domain $\Omega_+$. Then
(\ref{4.8.1}) implies that %
$$ 
|N_P^--N_Z^-|=|N_P^+-N_Z^+|\not=0.
$$ 

Below we include three  results on Blaschke products which follow
from our previous considerations. These results are familiar to
experts teaching Complex Analysis  and could be good exercises for
their students. The parts (1) and (2) will be used in the proof of
the converse statement of Theorem~3. The part (3) provides
additional information on critical points of related quadratic
differentials.

Let $A(z)=\prod_{k=1}^n\frac{z-a_k}{1-\overline{a}_k z}$ and
$B(z)=\prod_{k=1}^n\frac{z-b_k}{1-\overline{b}_k z}$ be Blaschke
products of order $n\ge 1$ each. Then the following holds.

\smallskip

\noindent %
\textbf{(1)} $A(z)$ (and $B(z)$) does not have critical points on
the unit circle $\mathbb{T}$; i.e. $A'(z)\not=0$ (and
$B'(z)\not=0$) if $|z|=1$. Therefore, $\mathbb{T}$ is a regular
trajectory of $Q_A(z)\,dz^2$ (and of $Q_B(z)\,dz^2$).

\smallskip

\noindent %
\textbf{(2)} $|\mathbb{T}|_{Q_A}=|\mathbb{T}|_{Q_B}=n$.  Therefore
$Q_A(z)\,dz^2$ and $Q_B(z)\,dz^2$ are coordinated in a sense of
Definition~2. %

\smallskip

\noindent %
\textbf{(3)} The quotient function $R(z)=\frac{A(z)}{B(z)}$ has a
critical point on $\mathbb{T}$; i.e. $R'(e^{i\theta})=0$ for some
$\theta$, $0\le \theta<2\pi$.

\medskip

In the proof of Theorem~3 below, we will use rational functions
depending on a variable $\zeta$.
Accordingly, we will use the following notations:  $R(\zeta)=C\frac{P_1(\zeta)}{P_2(\zeta)}$, where  $C\not=0$  and %
 \begin{equation} \label{4.10} %
P_1(\zeta)=\prod_{k=1}^{m_1}(\zeta-a_k)^{p_k}, \quad \quad
P_2(\zeta)=\prod_{k=1}^{m_2} (\zeta-b_k)^{q_k}.
 \end{equation} %
We assume in (\ref{4.10}) that $a_1,\ldots,a_{m_1}$,
$b_1,\ldots,b_{m_2}$ are distinct points on $\mathbb{C}$ and
$p_k$, $q_k$ are positive integers. Let $n_1=\sum_{k=1}^{m_1}
p_k$, $n_2=\sum_{k=1}^{m_2} q_k$. Calculating initial terms of the
power series expansion of $\frac{R'(\zeta)}{R(\zeta)}$ as in
(\ref{4.7.1}), we conclude that  $Q_R(\zeta)\,d\zeta^2$  has a
circular domain centered at $\zeta=\infty$ if and only if
$n_1\not=n_2$. Therefore, the domain configuration of
$Q_R(\zeta)\,d\zeta^2$ includes $N=m_1+m_2+1$ circle domains if
$n_1\not = n_2$ and it includes $N=m_1+m_2$ circle domains if
$n_1=n_2\ge 2$. In addition, the domain configuration of
$Q_R(\zeta)\,d\zeta^2$ may include a finite number of ring
domains.

Next, we show how to use quadratic differentials to prove
Theorem~3 for the rational function $R(\zeta)$ defined above. Then
Theorem~2 will follow as a special case.

\medskip

\noindent %
\emph{Proof of Theorem~3.}  Let
$R(\zeta)=C\frac{P_1(\zeta)}{P_2(\zeta)}$ be a rational function
of degree $n\ge 2$. Suppose further that $L_R(1)$ is analytic and
connected and that $n_1>n_2$. Then $n=n_1$ and $R(\infty)=\infty$.
Let $\Omega_-$ and $\Omega_+$ be inner and exterior domains for
$L_R(1)$. Then $\Omega_-$ contains all zeros $a_k$  of $R(\zeta)$
and $\Omega_+$ contains all poles $b_k$ of $R(\zeta)$ including
the pole at $\zeta=\infty$.

The lemniscate $L_R(1)$ is a regular closed trajectory of the
quadratic differential $Q_R(\zeta)\,d\zeta^2$ defined by
(\ref{4.5}). Transplanting $Q_R(\zeta)\,d\zeta^2$ from $\Omega_+$
and $\Omega_-$ to $\mathbb{D}_+$ and $\mathbb{D}$, respectively,
we obtain the quadratic differential %
\begin{equation}  \label{4.11} %
\begin{array}{rl} %
Q_{A_1}(z)\,dz^2&=Q_R(\varphi_+(z))(\varphi'_+(z))^2\,dz^2 \\ { } %
&
=-\frac{1}{4\pi^2}\left(\frac{R'(\varphi_+(z))}{R(\varphi_+(z))}\right)^2\left(\varphi'_+(z)\right)^2\,dz^2
 =
-\frac{1}{4\pi^2}\left(\frac{A'_1(z)}{A_1(z)}\right)^2\,dz^2,
\end{array} %
\end{equation} %
where $A_1(z)=R(\varphi_+(z))$ for $z\in \mathbb{D}_+$,
and the quadratic differential %
\begin{equation}  \label{4.11.1} %
\begin{array}{rl} %
Q_{B_1}(z)\,dz^2&=Q_R(\varphi_-(z))(\varphi'_-(z))^2\,dz^2 \\ {} %
&
=-\frac{1}{4\pi^2}\left(\frac{R'(\varphi_-(z))}{R(\varphi_-(z))}\right)^2\left(\varphi'_-(z)\right)^2\,dz^2
 = -\frac{1}{4\pi^2}\left(\frac{B'_1(z)}{B_1(z)}\right)^2\,dz^2,
\end{array} %
\end{equation} %
where $B_1(z)=R(\varphi_-(z))$ for $z\in \mathbb{D}$.

 Since
$|A_1(z)|=|R(\varphi_+(z))|=1$ for $|z|=1$ and since $A_1(z)$ does
not have zeroes in $\mathbb{D}_+$ it follows that $A_1(z)$ is a
Blaschke product of order $n$ with poles at the points
$\varphi_+^{-1}(b_k)\in \mathbb{D}_+$ and at $z=\infty$.
 Similarly, since $|B_1(z)|=|R(\varphi_-(z))|=1$ for $|z|=1$ and since $B_1(z)$ does
not have poles in $\mathbb{D}$ it follows that $B_1(z)$ is a
Blaschke product of order $n$ with zeros at the points
$\varphi_-^{-1}(a_k)\in \mathbb{D}$.

Applying equation~(\ref{2.4.3}) of Lemma~\ref{Lemma-1}, we
conclude that the fingerprint $k(z)$, $z=e^{i\theta}$,  can be
found as a solution to the
equation %
\begin{equation}  \label{4.11.2}
\log (A_1(k(z))/A_1(k(z_0)))=\log (B_1(z)/B_1(z_0)).%
\end{equation} %
 Since
$|A_1(k(z_0))|=1$ and $|B_1(z_0)|=1$, the latter equation implies
that $k(z)$ satisfies equation (\ref{1.4}) with Blaschke products
$A(z)=A_1(z)/A_1(k(z_0))$ and $B(z)=B_1(z)/B_1(z_0)$ of order $n$
each. This proves the first part of Theorem~3.

To prove the converse statement of Theorem~3, we consider two
Blaschke products $A(z)=\prod_{k=1}^n
\frac{z-\alpha_k}{1-\overline{\alpha}_k z}$ and
$B(z)=\prod_{k=1}^n \frac{z-\beta_k}{1-\overline{\beta}_k z}$ of
order $n\ge 2$ each such that $A(\infty)=\infty$. Let
$Q_A(z)\,dz^2$ and $Q_B(z)\,dz^2$ be quadratic differentials
defined by (\ref{4.5}). As we had mentioned above in the remark
(1), Blaschke products do not have critical points on the unit
circle. Thus, $\mathbb{T}$ is a regular trajectory of each of the
quadratic differentials $Q_A(z)\,dz^2$ and $Q_B(z)\,dz^2$.
Furthermore it follows from the remark (2) above that the
quadratic differentials $Q_A(z)\,dz^2$ and $Q_B(z)\,dz^2$ are
coordinated on $\mathbb{T}$.

Equation (\ref{1.4}) is equivalent to  equation (\ref{2.4.3}) with
${\mathcal{A}}$ and ${\mathcal{B}}$ defined by  (\ref{2.4.5}) with
$Q_+(z)\,dz^2=Q_A(z)\,dz^2$ and $Q_-(z)\,dz^2=Q_B(z)\,dz^2$. Thus,
if $k(z)$ is a solution to (\ref{1.4}) then $k(z)$ is a solution
to (\ref{2.4.3}). Now, since $Q_A(z)\,dz^2$ and $Q_B(z)\,dz^2$ are
coordinated on $\mathbb{T}$, it follows from
Theorem~\ref{Theorem-4}
 that there is a quadratic differential
$Q(\zeta)\,d\zeta^2$ defined on $\overline{\mathbb{C}}$ and its
regular trajectory $\Gamma$ such that %

\begin{equation}  \label{4.12} %
Q(\zeta)\,d\zeta^2= \left\{\begin{array}{ll} %
Q_{A}(\tau_+(\zeta))\left(\tau'_+(\zeta)\right)^2\,d\zeta^2 &\quad
{\mbox{for $\zeta\in \Omega_+\cup
\Gamma$,}}\\
Q_{B}(\tau_-(\zeta))\left(\tau'_-(\zeta)\right)^2\,d\zeta^2 &\quad
{\mbox{for $\zeta\in \Omega_-\cup \Gamma$,}} %
\end{array} %
\right. %
\end{equation} %
where $\tau_+:\Omega_+\to \mathbb{D}_+$ and $\tau_-:\Omega_-\to
\mathbb{D}$ are inverses of the conformal mappings $\varphi_+$ and
$\varphi_-$, respectively. Thus, the quadratic differentials
$Q_A(z)\,dz^2$ and $Q_B(z)\,dz^2$ are pullbacks of the quadratic
differential $Q(\zeta)\,d\zeta^2$ and therefore, by
Theorem~\ref{Theorem-3.1.1}, the fingerprint
$k=\varphi_+^{-1}\circ \varphi_-$ of $\Gamma$  given by equation
(\ref{2.4.3}) or, equivalently, by equation (\ref{1.4}).

It still remains to show that the trajectory $\Gamma$ defined
above is a lemniscate of a rational function with required
properties.
Equations (\ref{4.12}) and (\ref{4.11}), (\ref{4.11.1}) imply that %
 \begin{equation}  \label{4.13} %
Q(\zeta)\,d\zeta^2= \left\{\begin{array}{ll} %
-\frac{1}{4\pi^2}
\left(\frac{A'_1(\zeta)}{A_1(\zeta)}\right)^2\,d\zeta^2 &\quad
{\mbox{for $\zeta\in \Omega_+\cup
\Gamma$,}}\\
-\frac{1}{4\pi^2}
\left(\frac{B'_1(\zeta)}{B_1(\zeta)}\right)^2\,d\zeta^2 &\quad
{\mbox{for $\zeta\in \Omega_-\cup \Gamma$,}} %
\end{array} %
\right. %
 \end{equation} %
where $A_1(\zeta)=A(\tau_+(\zeta))$ and
$B_1(\zeta)=B(\tau_-(\zeta))$.

Since $Q(\zeta)\,d\zeta^2$ is a quadratic differential on
$\overline{\mathbb{C}}$, $Q(\zeta)$ is a rational function. It
follows from (\ref{4.13}) that all zeros and poles of $Q(\zeta)$
are of even orders and therefore there is a rational function
$R_1(\zeta)$ such that %
$$ %
Q(\zeta)\,d\zeta^2=-\frac{1}{4\pi^2}R_1^2(\zeta)\,d\zeta^2. %
$$ %
Furthermore, since
$R_1^2(\zeta)=\left(A'_1(\zeta)/A_1(\zeta)\right)^2$ for $\zeta\in\Omega_+$  
and $R_1^2(\zeta)=\left(B'_1(\zeta)/B_1(\zeta)\right)^2$ for $\zeta\in \Omega_-$ 
it follows that all finite poles of $R_1(\zeta)$ are simple with
real integer residues. Also, $A_1(\zeta)=A(\tau_+(\zeta))$ has a
pole at $\zeta=\infty$ and therefore $\zeta=\infty$ is a simple
zero of $R_1(\zeta)$ and the limit $\lim_{\zeta\to \infty} \zeta
R_1(\zeta)$ is a positive integer. The latter information implies
that the
function $R(\zeta)$ defined as%
\begin{equation}  \label{4.13.1} %
R(\zeta)=\exp\left(\int_{\zeta_0}^\zeta R_1(t)\,dt\right)
 \end{equation} %
is a rational function such that $R'(\zeta)/R(\zeta)=R_1(\zeta)$
and $R(\infty)=\infty$.
This shows that %
\begin{equation} \label{4.14} %
Q(\zeta)\,d\zeta^2=Q_R(\zeta)\,d\zeta^2=-\frac{1}{4\pi^2}\left(\frac{R'(\zeta)}{R(\zeta)}\right)^2\,d\zeta^2
\end{equation}  %
and therefore $\Gamma$ is a lemniscate of the rational function
$R(\zeta)$. In particular, if in (\ref{4.13.1}) we choose
$\zeta_0\in \Gamma$ then $R(\zeta_0)=1$ and therefore
$\Gamma=L_R(1)$. This completes the proof of Theorem~3.  \hfill
$\Box$





\medskip

\noindent %
\emph{Proof of Theorem~2.} Taking $P_2(\zeta)\equiv 1$ in the
first part of the proof of Theorem~3 given above, we conclude that
(\ref{4.11.2}) holds with $A_1(z)=c z^n$ where $|c|=1$.  Thus,
(\ref{4.11.2}) is equivalent to (\ref{1.2}) in this case.

Similarly, taking $A(z)=z^n$ in the second part of the proof of
Theorem~3 given above, we conclude that (\ref{4.14}) holds with
the rational function $R(\zeta)$ whose the only pole is at
$\infty$. Hence, $R(\zeta)$ is a polynomial of order $n$ in this
case.

The uniqueness statement of the converse part of Theorem~2 easily
follows from the maximum modulus principle. Indeed, if there are
two polynomials $P_1(\zeta)$ and $P_2(\zeta)$ of degree $n\ge2$
with positive leading coefficients, which share the lemniscate
$\Gamma$, then their quotient $P_1(\zeta)/P_2(\zeta)$ is analytic
and non-vanishing  on $\Omega_+$ (including the point
$\zeta=\infty$) and satisfies the equation
$|P_1(\zeta)/P_2(\zeta)|=1$ on $\Gamma=\partial \Omega_+$. Hence,
$P_1(\zeta)=P_2(\zeta)$ as required. \hfill $\Box$

\medskip

We note here that equation~(\ref{4.13.1}) defines a unique
rational function $R(\zeta)$ having $\Gamma$ as its lemniscate
$L_R(1)$. However, (unlike in Theorem~\ref{Theorem-2}) a rational
function $R(\zeta)$, which existence is quarantined by the
converse part of Theorem~\ref{Theorem-3}, is not unique in
general. Indeed, any two Blaschke products  $B_1(\zeta)$ and
$B_2(\zeta)$ of degree $n$ such that $B_1(\infty)=\infty$,
$B_2(\infty)=\infty$ share the unit circle $\mathbb{T}$  as their
lemniscates $L_{B_1}(1)$ and $L_{B_2}(1)$. Whether or not there
are non-circular rational lemniscates shared by two essentially
distinct rational functions of the same degree remains an open
question to this author.

\section{Polygonal curves and related quadratic differentials}  %
\setcounter{equation}{0}  %

In this section, we discuss two particular cases of quadratic
differentials $Q_-(z)\,dz^2$ and $Q_+(z)\,dz^2$ and corresponding
curves $\Gamma$ generated by fingerprints $k(z)$ defined by
formulas (\ref{2.4.3}) and (\ref{2.4.5}).

\medskip

\textbf{(a) Cartesian polygonal curves.} By a Cartesian polygonal
curve we  understand a Jordan curve consisting of a finite number
of horizontal and vertical segments. Any such curve $\Gamma$ is a
boundary of a \emph{standard polygon} $\Omega_-$ having an even
number of sides and even number of vertices, $v_1,\ldots,v_{2n}$.
We suppose here that vertices are always oriented in the
counterclockwise direction and that $v_{2n+1}=v_1$, $v_0=v_{2n}$.
We assume additionally that the boundary segment $[v_1,v_2]$
represents a horizontal side of $\Omega_-$. Let $\alpha_k\pi$ be
the angle of $\Omega_-$ at its vertex $v_k$.
 Then either
$\alpha_k=\frac{1}{2}$ or $\alpha_k=\frac{3}{2}$. Since the sum of
angles of any polygon with $2n$ vertices  is $2\pi(n-1)$, one can
easily see that $\Omega_-$ must have $n+2$ vertices with angles
$\frac{\pi}{2}$ and $n-2$ vertices with angles $\frac{3\pi}{2}$.

The horizontal and vertical sides of $\Omega_-$ are arcs of
trajectories and, respectively,  arcs of orthogonal trajectories
of the quadratic differential $Q(\zeta)\,d\zeta^2=1\cdot
d\zeta^2$. Transplanting this quadratic differential via the
mapping $\varphi_-:\mathbb{D}\to \Omega_-$,  we obtain the
following quadratic  differential:
\begin{equation} \label{5.1} %
Q_-(z)\,dz^2=  C_- e^{i\gamma_-}\prod_{k=1}^{2n}
(z-e^{i\beta_k^-})^{2(\alpha_k-1)} \,dz^2, \quad \quad z\in
\mathbb{D},
\end{equation} %
with some  $C_->0$, $\gamma_-\in \mathbb{R}$, and with
$e^{i\beta_k^-}=\tau_-(v_k)$, where  $0\le
\beta_1^-<\beta_2^-<\cdots<\beta_{2n}^-<\beta_1^-+2\pi$. Figure~6
presents an example of a Cartesian polygonal curve $\Gamma$ and
shows critical trajectories of the corresponding quadratic
differential $Q_-(z)\, dz^2$.

Since the unit circle $\mathbb{T}$ consists of arcs of
trajectories and/or arcs of orthogonal trajectories of
$Q_-(z)\,dz^2$, the quadratic differential (\ref{5.1}) is real in
all its regular points of $\mathbb{T}$. This after simple
calculation
gives the following value for $\gamma_-$: %
$$ %
\gamma_-=\sum_{k=1}^{2n} (1-\alpha_k)\beta_k^-.  %
$$ %

Now, the value of a positive constant $C_-$, which does not affect
the trajectory structure of $Q_-(z)\,dz^2$, can be obtained from
the equation for the $Q_-$-length of the arc
$\{e^{i\theta}:\,\beta_1^-\le
\theta\le \beta_2^-\}$, which gives the following: %
$$ %
C_-=-e^{-i\gamma_-}|v_2-v_1|^2 \left(\int_{\beta_1^-}^{\beta_2^-}
\prod_{k=1}^{2n}(e^{i\theta}-e^{i\beta_k^-})^{\alpha_k-1}e^{i\theta}\,d\theta\right)^{-2}. %
$$ %

Similarly, transplanting $Q(\zeta)\,d\zeta^2=1\cdot d\zeta^2$ via
the  mapping  $\varphi_+:\mathbb{D}_+\to \Omega_+$, we obtain the
second required quadratic differential:%
\begin{equation} \label{5.2} %
Q_+(z)\,dz^2=  C_+e^{i\gamma_+}z^{-4}\prod_{k=1}^{2n}
(z-e^{i\beta_k^+})^{2(1-\alpha_k)} \,dz^2, \quad \quad z\in
\mathbb{D}_+,
\end{equation} %
where   $e^{i\beta_k^+}=\tau_+(v_k)$ with $0\le
\beta_1^+<\beta_2^+<\cdots<\beta_{2n}^+<\beta_1^+ +2\pi$. The
constants $\gamma_+\in \mathbb{R}$ and $C_+>0$ in (\ref{5.2}) are
given by equations: %
$$ %
\gamma_+=\sum_{k=1}^{2n} (\alpha_k-1)\beta_k^+
$$ %
and %
$$ %
C_+=-e^{-i\gamma_+}|v_2-v_1|^2 \left(\int_{\beta_1^+}^{\beta_2^+}
\prod_{k=1}^{2n}(e^{i\theta}-e^{i\beta_k^+})^{1-\alpha_k}e^{-i\theta}\,d\theta\right)^{-2}. %
$$ %

We note that quadratic differentials (\ref{5.1}) and (\ref{5.2})
are coordinated in the sense of Definition~2 if and only if the
following equations are satisfied for all $k=1,\ldots,2n$: %
\begin{equation} \label{5.3} %
\frac{\int_{\beta_{k-1}^+}^{\beta_k^+}\prod_{j=1}^{2n}\left(e^{i\theta}-e^{i\beta_j^+}\right)^{1-\alpha_j}\,e^{-i\theta}\,d\theta}
{\int_{\beta_{k-1}^-}^{\beta_k^-}\prod_{j=1}^{2n}\left(e^{i\theta}-e^{i\beta_j^-}\right)^{\alpha_j-1}\,e^{i\theta}\,d\theta}=Ce^{i\gamma},%
\end{equation} %
where $C=\sqrt{C_-/C_+}$ and $\gamma=(\gamma^--\gamma^+)/2$.

Finally, if $Q_-(z)\,dz^2$ and $Q_+(z)\,dz^2$ given by formulas
(\ref{5.1}), (\ref{5.2}) are coordinated, i.e. if they satisfy
equation (\ref{5.3}), then by Theorem~6 they define a polygonal
curve $\Gamma$ whose fingerprint $k:\mathbb{T}\to \mathbb{T}$  can
be found from equations (\ref{2.4.3}) and  (\ref{2.4.5}) for the
quadratic differentials $Q_-(z)\,dz^2$ and $Q_+(z)\,dz^2$ given by
 (\ref{5.1}) and (\ref{5.2}), respectively.


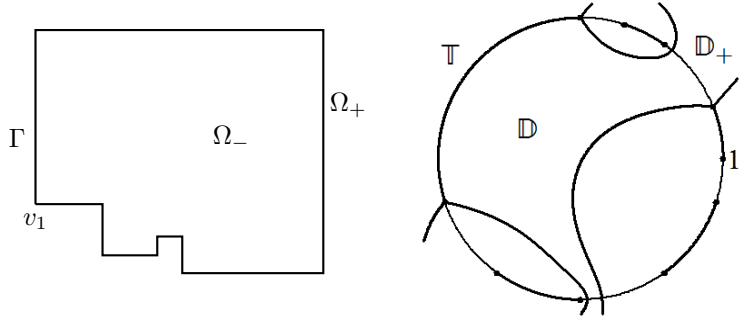
\begin{figure}
\hspace{-3.5cm} %
\begin{minipage}{0.10\linewidth}
\begin{tikzpicture} [thick,scale=0.22] 
\draw [black, thick] (0,0) --(4.00956,0) -- (4.00956,-3.08468) --
(4.00956+3.27259,-3.08468)--(4.00956+3.27259,-3.08468+1.13451)
--(4.00956+3.27259+1.49598,-3.08468+1.13451)--(4.00956+3.27259+1.49598,-3.08468+1.13451-2.20908)--
(4.00956+3.27259+1.49598+8.43538,-3.08468+1.13451-2.20908)--(4.00956+3.27259+1.49598+8.43538,-3.08468+1.13451-2.20908+14.6512)--
(4.00956+3.27259+1.49598+8.43538-17.2135,-3.08468+1.13451-2.20908+14.6512)--
(4.00956+3.27259+1.49598+8.43538-17.2135,-3.08468+1.13451-2.20908+14.6512-10.492);

\node at (0,4) [left] {$\Gamma$};  %
\node at (0,0) [below] {$v_1$}; %
 \node at (10,4) [right] {$\Omega_-$};
 \node at (17,6) [right] {$\Omega_+$};

\end{tikzpicture} %
\end{minipage} %
\hspace{4cm} %
\begin{minipage}{0.10\linewidth}
$$\includegraphics[scale=.45,angle=0]{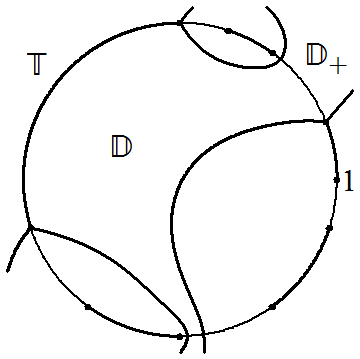}
$$
\end{minipage} %
\vspace{0cm} %
\caption{Cartesian polygonal curve and critical trajectories of $Q_-(z)\,dz^2$.} %
\end{figure}



\FloatBarrier

\medskip %

\textbf{(b) Polar polygonal curves.} We start with the quadratic
differential %
\begin{equation} \label{5.4.1} %
 Q(\zeta)\,d\zeta^2=-\frac{d\zeta^2}{\zeta^2}. %
\end{equation} %
Then the radial  segments of the form
$\{\zeta=re^{i\alpha}:\,r_1\le r\le r_2\}$ with some $\alpha\in
\mathbb{R}$ and $0<r_1<r_2<\infty$ are closed arcs on the
orthogonal trajectories of $Q(\zeta)\,d\zeta^2$ and the closed
arcs of circles centered at $\zeta=0$ are closed arcs on the
trajectories of $Q(\zeta)\,d\zeta^2$.

By a polar polygonal curve $\Gamma$ we mean a closed Jordan curve
bounded by a finite number of radial segments and circular arcs as
above. The components $\Omega_-$ and $\Omega_+$ of
$\overline{\mathbb{C}}\setminus \Gamma$ will be called a polar
polygon and an outer polar polygon, respectively. We will assume
additionally that $0\in \Omega_-$. We want to mention here that
polar polygons, also known as \emph{gearlike domains}, were used
by R.W.~Barnard and his collaborators in their study of Goodman's
omitted area problem; see, for instance,~\cite{BP}.

 Each
of the polygons $\Omega_-$ and $\Omega_+$ has an even number of
sides, say $2n$, and accordingly they have an even number of
vertices and same even number of corresponding angles. For all
these objects we will use notations introduced in part
\textbf{(a)} of this section. Since $0\in \Omega_-$ one can easily
find that $\Omega_-$ has $n$ vertices $v_k$ with angles
$\alpha_k\pi=\frac{1}{2}\pi$ and $n$ vertices $v_k$ with angles
$\alpha_k\pi=\frac{3}{2}\pi$.

Transplanting $Q(\zeta)\,d\zeta^2$ via the mapping
 $\varphi_-:\,\mathbb{D}\to \Omega_-$  and assuming that $\varphi(0)=0$, we obtain the
 following quadratic differential:
\begin{equation} \label{5.5} %
Q_-(z)\,dz^2=  -C_- e^{i\gamma_-}z^{-2} \prod_{k=1}^{2n}
(z-e^{i\beta_k^-})^{2(\alpha_k-1)} \,dz^2, \quad \quad z\in
\mathbb{D},
\end{equation} %
where  $e^{i\beta_k^-}=\tau_-(v_k)$ with $0\le
\beta_1^-<\beta_2^-<\cdots<\beta_{2n}^-<\beta_1^-+2\pi$. Figure~7
displays  an example of a polar polygonal curve $\Gamma$ and shows
critical trajectories of the corresponding quadratic differential
$Q_-(z)\, dz^2$.

Comparing the coefficients for term $\zeta^{-2}$ in (\ref{5.4.1})
and for the term  $z^{-2}$ in (\ref{5.5}), we obtain the following
equation: %
$$ %
C_-e^{i\left(\gamma_- +
2\sum_{k=1}^n(\alpha_k-1)\beta_k^-\right)}=1.
$$ %
The latter implies that $C_-=1$ and %
$$ %
\gamma_-=2\sum_{k=1}^{2n} (1-\alpha_k)\beta_k^-.
$$ %

Next, the unit circle $\mathbb{T}$ consists of arcs of
trajectories and orthogonal trajectories of $Q_-(z)\,dz^2$. Hence
$Q_-(z)\,dz^2$ is real at all its regular points on $\mathbb{T}$.
This, after some algebra,  leads to the following condition: %
$$ %
e^{i\gamma_-}=-e^{2i\sum_{k=1}^{2n} (1-\alpha_k)\beta_k^-}.
$$ %
Furthermore, we will assume that the boundary arc of $\partial
\Omega_-$ with endpoints $v_1$ and $v_2$ is a circular arc. Then
the quadratic differential $Q_-(z)\,dz^2$ must be positive on the
arc $\{e^{i\theta}:\,\beta_1^-<\theta<\beta_2^-\}$. The latter
condition, after routine calculation, implies the following: %
\begin{equation}  \label{5.8} %
\sum_{k=1}^{2n}(1-\alpha_k)\beta_k^-=\pi \quad \quad ({\mbox{mod
$2\pi$}}).
\end{equation}%
This equation is equivalent to the following %
\begin{equation} \label{5.9} %
{\sum}_1\beta_k^- -{\sum}_2 \beta_k^-=2\pi \quad \quad ({\mbox{mod
$4\pi$}}),
\end{equation} 
where the first sum is taken over all $k$ such that
$e^{i\beta_k^-}$ is a zero of $Q_-(z)\,dz^2$ and the second sum is
taken over all $k$ such that $e^{i\beta_k^-}$ is a pole of
$Q_-(z)\,dz^2$.

Combining our observations, we conclude that $Q_-(z)\,dz^2$ can be
represented in the form %
$$ %
Q_-(z)\,dz^2=  -z^{-2} \prod_{k=1}^{2n}
(z-e^{i\beta_k^-})^{2(\alpha_k-1)} \,dz^2
$$ %
with $\beta_1^-,\ldots,\beta_{2n}^-$ satisfying equation
(\ref{5.8}) or, equivalently, equation (\ref{5.9}).

The same argument as above shows that the quadratic differential
$Q_+(z)\,dz^2$ has the form %

 \begin{equation} \label{5.6} %
Q_+(z)\,dz^2=  -z^{-2}\prod_{k=1}^{2n}
(z-e^{i\beta_k^+})^{2(1-\alpha_k)} \,dz^2, \quad \quad z\in
\mathbb{D}_+,
\end{equation} %
where   $e^{i\beta_k^+}=\tau_+(v_k)$ with  $0\le
\beta_1^+<\beta_2^+<\cdots<\beta_{2n}^+<\beta_1^+ +2\pi$.
Furthermore, the parameters $\beta_1^+,\ldots,\beta_{2n}^+$
satisfy the following equation: %
$$ %
\sum_{k=1}^{2n}(\alpha_k-1)\beta_k^+=\pi \quad \quad ({\mbox{mod
$2\pi$}}).
$$ %

As in the previous case,  quadratic differentials (\ref{5.5}) and
(\ref{5.6}) are coordinated in the sense of Definition~2 if and
only if  %
\begin{equation} \label{5.7} %
\frac{\int_{\beta_{k-1}^+}^{\beta_k^+}\prod_{j=1}^{2n}\left(e^{i\theta}-e^{i\beta_j^+}\right)^{1-\alpha_j}\,d\theta}
{\int_{\beta_{k-1}^-}^{\beta_k^-}\prod_{j=1}^{2n}\left(e^{i\theta}-e^{i\beta_j^-}\right)^{\alpha_j-1}\,d\theta}=1
\quad \quad  {\mbox{for $k=1,2, \ldots, 2n$.}}%
\end{equation} %

As in part \textbf{(a)}, we conclude that if the quadratic
differentials  $Q_-(z)\,dz^2$ and $Q_+(z)\,dz^2$ given by formulas
(\ref{5.5}) and  (\ref{5.6}) are coordinated, i.e. if they satisfy
equations (\ref{5.7}), then  they define a polar polygonal curve
$\Gamma$ whose fingerprint $k:\mathbb{T}\to \mathbb{T}$  can be
found from equations (\ref{2.4.3}) and (\ref{2.4.5}).

\FloatBarrier
\begin{figure} %
\hspace{-0.5cm} %
\begin{minipage}{.5\textwidth}
\begin{tikzpicture} [thick,scale=1.2] , 
\draw [black, thick] (0.18395320871404977,
0.9829350014135244)--(0.3610714038274868, 1.9293478124823114);

\draw [black, thick] (-0.09740313413217537,
-1.960425507295154)--(-0.024140230259961204,
-0.48586858703526326);

\draw [black, thick] (-0.13502406259764815,
-0.46735376012638813)--(-0.1187557092708051,
-0.41104471452300717);

\draw [black, thick] (0.0815743,-0.420008)--(0.1400651133227818,
-0.7211631479541237);

\draw [black, thick] (0.7346390419363942, 0)--(1,0);

 \draw [black,thick,domain=0:1.3857895357235062*180/pi] plot ({cos(\x)}, {sin(\x)});

  \draw [black,thick,domain=1.3857895357235062*180/pi:(1.3857895357235062 + 3.2769555772941956)*180/pi] plot ({e^(0.674394308508077)*cos(\x)}, {e^(0.674394308508077)*sin(\x)});

  \draw [black,thick,domain=(1.3857895357235062 + 3.2769555772941956)*180/pi:(1.3857895357235062 + 3.2769555772941956 - 0.23160957260640253)*180/pi]
  plot ({e^(0.674394308508077- 1.3949786339316466)*cos(\x)}, {e^(0.674394308508077- 1.3949786339316466)*sin(\x)});

 \draw [black,thick,domain=((1.3857895357235062 + 3.2769555772941956 - 0.23160957260640253)*180/pi:(1.3857895357235062 + 3.2769555772941956 - 0.23160957260640253+0.47308625844850927)*180/pi]
 plot ({e^(0.674394308508077- 1.3949786339316466-0.12838448419039247)*cos(\x)}, {e^(0.674394308508077- 1.3949786339316466-0.12838448419039247)*sin(\x)});

\draw [black,thick,domain=((1.3857895357235062 +
3.2769555772941956 -
0.23160957260640253+0.47308625844850927)*180/pi:360]
 plot ({e^(0.674394308508077- 1.3949786339316466-0.12838448419039247+0.5405928098095882)*cos(\x)}, {e^(0.674394308508077- 1.3949786339316466-0.12838448419039247+0.5405928098095882)*sin(\x)});

\node at (0.92,0.95) [left] {$\Gamma$};  %
\node at (1,0) [right] {$v_1$}; %
\node at (0,0) [right] {$0$}; %
\node at (0,0)  {$\bullet$}; %
 \node at (-1,0) [right] {$\Omega_-$};
 \node at (0.35,-1.2) [right] {$\Omega_+$};



\end{tikzpicture} %
\end{minipage} %
\hspace{-2cm} %
\begin{minipage}{.5\textwidth}

$$\includegraphics[scale=.46,angle=0]{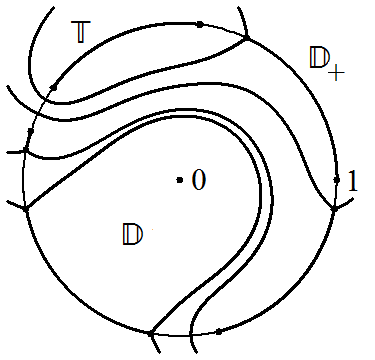}
$$
\vspace{0cm} %
\end{minipage} %
\caption{Polar polygonal curve and critical trajectories of $Q_-(z)\,dz^2$.} %
\end{figure}
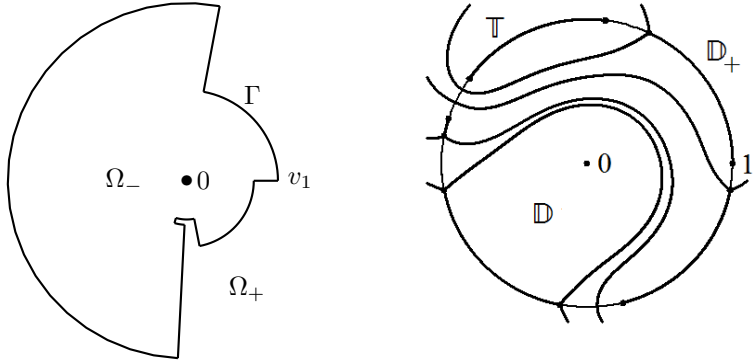


\medskip

\textbf{(c)}  A few remarks are in order now.

(1) The mapping functions $\varphi_-$ and $\varphi_+$ used in
parts \textbf{(a)} and \textbf{(b)} of this section can be
represented by the well-known Schwarz-Christoffel integrals. In
fact, our formulas (\ref{5.3}) and (\ref{5.7}) can be obtained as
Schwarz-Christoffel integrals taken over the corresponding arcs of
the unit circle.

(2) A simple fact of geometry that, in the case \textbf{(a)}, the
polygon $\Omega_-$ has $n+2$ vertices with angle $\frac{\pi}{2}$
and $n-2$ vertices with angles $\frac{3\pi}{2}$ can be interpreted
as an elementary corollary of the following formula: %
$$ %
p-q=4-4g-2s,
$$ %
see, for example,  Lemma~3.2 in \cite{J2}). This formula relates
the number of poles $p$ and number of zeros $q$ (both counted with
multiplicities) of the quadratic differential $Q(z)\,dz^2$ defined
on a domain $D$ lying on a compact Riemann surface of genus $g$
when $\partial D$ consists of $s$ non-degenerate curves and
$Q(z)\,dz^2>0$ on $\partial D$. Of course, in the case of the
complex sphere this formula reduces to $p-q=4$.

Now, since the vertices $v_k$ of $\Omega_-$ with angles
$\frac{\pi}{2}$ correspond to simple poles of $Q_-(z)\,dz^2$ and
the vertices $v_k$ with angles $\frac{3\pi}{2}$ correspond to
simple zeroes of $Q_-(z)\,dz^2$, we conclude from (\ref{5.1}) that
the difference between the number of vertices of $\Omega_-$ with
angles $\frac{\pi}{2}$ and the number of its vertices with angles
$\frac{3\pi}{2}$ equals $4$. The same conclusion can be made by
comparing the number of poles and zeroes of the quadratic
differential (\ref{5.2}).

The same argument, being applied to the quadratic differentials
(\ref{5.5}) and (\ref{5.6}), shows that, in the case of polar
polygons, the number of vertices of $\Omega_-$ with angle
$\frac{\pi}{2}$ should be equal to the number of its vertices with
angle $\frac{3\pi}{2}$.

(3) Equations (\ref{5.3}) and (\ref{5.7}) give necessary and
sufficient conditions which guarantee that the Schwarz-Christoffel
integrals representing functions $\varphi_-$ and $\varphi_+$
 define one-to-one mappings from $\mathbb{D}$ and
$\mathbb{D}_+$ onto polygons $\Omega_-$ and $\Omega_+$,
respectively. Of course, experts know that a similar fact holds
true for the Schwarz-Christoffel mappings from $\mathbb{D}$ and
$\mathbb{D}_+$ onto any two complementary polygons with common
Jordan boundary.

Surprisingly to  this author, the latter fact is not mentioned in
standard  textbooks on Complex Analysis. Thus, we  state it here.

\medskip

\noindent %
\begin{proposition}For $n\ge 3$, let $0\le \beta_1^-<\beta_1^-<\cdots
<\beta_n^-<\beta_1^-+2\pi$ and let $0<\alpha_k<2$,
$k=1,2,\ldots,n$, be such that $\sum_{k=1}^n \alpha_k=n-2$.

Then the Schwarz-Christoffel integral %
$$ %
F(z)=\int_0^z \prod_{k=1}^n (\tau-e^{i\beta_k^-})^{\alpha_k-1}\,d\tau %
$$ %
maps $\mathbb{D}$ conformally and one-to-one onto some polygon if
and only if there are points $z_k^+=e^{i\beta_k^+}$ with
$0\le\beta_1^+<\beta_1^+<\cdots <\beta_n^+<\beta_1^++2\pi$ such
that the equations (\ref{5.3}) with some $C>0$ and $\gamma\in
\mathbb{R}$ are satisfied for all $k=1,2,\ldots,n$.
\end{proposition}

\end{document}